\documentclass[twoside]{article}

\usepackage[accepted]{aistats2025}
\usepackage{amssymb}
\usepackage{float}          
\usepackage{mathtools}
\usepackage{amsmath}
\usepackage{amsthm}
\usepackage[dvipsnames]{xcolor}
\usepackage{multirow}
\usepackage{multicol}
\usepackage{booktabs}
\usepackage{diagbox}
\usepackage{varioref}
\usepackage[unicode=true,psdextra,colorlinks]{hyperref}
\hypersetup{%
      breaklinks,
      hypertexnames = true,
      plainpages    = false,
      colorlinks = true,
      urlcolor   = BrickRed,
      linkcolor  = MidnightBlue,
      filecolor  = Mulberry,      
      citecolor  = BrickRed,
      menubordercolor = BrickRed
    }
\usepackage[noabbrev, capitalise, nameinlink]{cleveref}
\usepackage{ellipsis}

\newtheorem{theorem}{Theorem}
\newtheorem{proposition}{Proposition}
\newtheorem{example}{Example}
\newtheorem{corollary}{Corollary}
\newtheorem{conjecture}{Conjecture}
\newtheorem{remark}{Remark}
\newtheorem{definition}{Definition}
\newtheorem{assumption}{Assumption}
\newtheorem{lemma}{Lemma}

\Crefname{conjecture}{Conjecture}{Conjectures}
\Crefname{assumption}{Assumption}{Assumptions}

\usepackage{enumitem}
\setlist[itemize]{align=parleft,left=0pt..1em}

\usepackage{makecell}

\usepackage{pifont}
\newcommand{\cmark}{\ding{51}}%
\newcommand{\xmark}{\ding{55}}%
\usepackage{bm}
\usepackage{ascii}
\newcommand{\centered}[1]{\begin{tabular}{l} #1 \end{tabular}}
\usepackage{comment}

\DeclareMathOperator*{\argmin}{argmin}
\DeclareMathOperator*{\argmax}{argmax}
\DeclareMathOperator*{\maximize}{maximize}
\DeclareMathOperator*{\minimize}{minimize}
\DeclareMathOperator{\prox}{prox}

\DeclareMathOperator*{\sgn}{sgn}
\DeclareMathOperator*{\dom}{dom}
\DeclareMathOperator*{\ri}{ri}
\DeclareMathOperator*{\range}{range}

\usepackage{tikz,forest}
  \usetikzlibrary{shapes,arrows,fit,calc,positioning} %
  \tikzset{box/.style={draw, rectangle, thick, text centered, minimum height=0.5cm, minimum width=1cm}}
  \tikzset{line/.style={draw, thick, -latex'}}

\usepackage{xparse}

\usepackage[round]{natbib}

\begin{document}

%
\runningtitle{Tight analysis of DCA improves convergence rates for Proximal Gradient Descent}

%

\twocolumn[

\aistatstitle{Tight Analysis of Difference-of-Convex Algorithm (DCA) Improves Convergence Rates for Proximal Gradient Descent}

\aistatsauthor{ Teodor Rotaru \And Panagiotis Patrinos \And  Fran\c{c}ois Glineur }

\aistatsaddress{ KU Leuven, UCLouvain \And  KU Leuven \And UCLouvain } ]

\begin{abstract}
We investigate a difference-of-convex (DC) formulation where the second term is allowed to be weakly convex. We examine the precise behavior of a single iteration of the difference-of-convex algorithm (DCA), providing a tight characterization of the objective function decrease, distinguishing between six distinct parameter regimes. %
Our proofs, inspired by the performance estimation framework, are notably simplified compared to related prior research. We subsequently derive sublinear convergence rates for the DCA towards critical points, assuming at least one of the functions is smooth. %
Additionally, we explore the underexamined equivalence between proximal gradient descent (PGD) and DCA iterations, demonstrating how DCA, a parameter-free algorithm, without the need for a stepsize, serves as a tool for studying the exact convergence rates of PGD. 
Finally, we propose a method to optimize the DC decomposition to achieve optimal convergence rates, potentially transforming the subtracted function to become weakly convex.
\end{abstract}

\section{\MakeUppercase{Introduction}}
Consider the difference-of-convex formulation
\begin{align}\label{eq:DC_program}\tag{DC}
    \minimize_{x \in \mathbb{R}^d} F(x)  {}\coloneqq{}  f_1(x) - f_2(x),
\end{align}
where $f_1, f_2: \mathbb{R}^d \rightarrow \mathbb{R}$ are proper, lower semicontinuous convex functions, and $F$ is lower bounded. 

A standard method to solve \eqref{eq:DC_program} is the difference-of-convex algorithm (DCA), a versatile method with no parameter that can find a \emph{critical} point of $F$, defined as a point $x^* \in \mathbb{R}^d$ for which there exists subgradients $g_1^* \in \partial 
f_1(x^*)$ and $g_2^* \in \partial 
f_2(x^*)$ such that $g_1^*=g_2^*$. 
Stationary points of $F$ are always critical, but the converse is not true.
Extensive analyses of DCA are provided by Dinh and Thi \citeyearpar{Dinh_Thi_dca_1997}, Tao and An \citeyearpar{DCA_Trust_region_1998}, Horst and
Thoai \citeyearpar{Horst_Thoai_DC_overview_1999},  Le Thi and Pham Dinh \citeyearpar{LeThi_2018_30_years_dev}. DCA is also referred to as the convex-concave procedure (CCCP), as seen in the work of Yuille and Rangarajan \citeyearpar{CCCP_2001_init_Alan_Anand}, Lanckriet and Sriperumbudur \citeyearpar{CCCP_2009_Convergence},  Lipp
and Boyd \citeyearpar{CCCP_Lipp_Boyd_2016}. Interestingly, convergence analysis of many methods can be reduced to the one of DCA; for example, the Frank-Wolfe algorithm (Yurtsever and Sra \citeyearpar{Yurtsever_Suvrit_FW_CCCP_2022}) or the proximal gradient descent (PGD) (Le Thi and Pham Dinh \citeyearpar[Section 3.3.4]{LeThi_2018_30_years_dev}). Conversely, Faust et al. \citeyearpar{bregman_DCA_2023} show that DCA is an instance of the Bregman proximal point algorithm. 

An extensive list of DCA applications is provided by Le Thi and Pham Dinh \citeyearpar{LeThi_2018_30_years_dev}. Notable examples include efficient formulations for clustering problems (Hoai An et al. \citeyearpar{A1_HOAIAN2014388}), dictionary learning (Vo et al. \citeyearpar{A2}), robust support vector regression (Wang et al. \citeyearpar{A3_Wang2015}), multi-class support vector machines (MSVM) (Le Thi and Nguyen \citeyearpar{A4_LeThi2017}), sparse logistic regression (Yang and Qian \citeyearpar{A5_Yang2016}), compressed sensing (Yin
et al. \citeyearpar{A6_Lou_2015_l1_l2_compressed_sensing}), adversarial attack for adversarial robustness and approximate optimization of complex functions (Awasthi et al. \citeyearpar{A7_Awasthi2024}) or Shallow Multilayer Perceptron (MLP) Neural Networks (Askarizadeh et al.\citeyearpar{A8_AMTNK24_CC_Shallow_NN}). Sun et al. \citeyearpar{SLN_NN_Shortcuts_2024_DCA} introduce the Negative ResNets, where $f_2$ is weakly convex.

We derive convergence rates to critical points when at least one of $f_1$ and $f_2$ is smooth (namely continuously differentiable, with Lipschitz gradient). Abbaszadehpeivasti et al. \citeyearpar{abbaszadehpeivasti2021_DCA} provide exact convergence rates of DCA when both functions are convex. Their approach, based on performance estimation (PEP) introduced by Drori and Teboulle \citeyearpar{drori_performance_2014} and refined by Taylor et al. \citeyearpar{taylor_smooth_2017}, leads to rigorous proofs for those rates. Their exactness is supported by strong numerical evidence and, in some cases, by the identification of instances matching those rates exactly. 

In this work, we generalize the standard \eqref{eq:DC_program} setting and consider the case where $f_2$ can be weakly convex (or hypoconvex). Some previous works also introduce weak convexity in either $f_1$ (Sun and Sun \citeyearpar{DC_alg_dif_Moreau_Smoothing_2022}) or $f_2$ (Syrtseva et al. \citeyearpar{dc_weakly_cvx_bundle_method_2023_Syrtseva}). 

A key motivation for examining the case with $f_2$ weakly convex is that it mirrors the behavior of applying PGD with stepsizes larger than the inverse Lipschitz constant (see \cref{sec:PGD_equivalence}). Additionally, our generalized DCA framework provides a useful tool for analyzing exact rates for PGD, with the benefit of handling one fewer parameter - DCA involves four curvature parameters compared to PGD's four curvature parameters plus the stepsize. Therefore, due to the equivalence of the iterations, it is more convenient to use a DCA-like analysis. 
%
%
Our results follow the same line as Abbaszadehpeivasti et al. \citeyearpar{abbaszadehpeivasti2021_DCA}, also relying on performance estimation. More precisely:%
\begin{itemize}
\item We characterize in \cref{thm:one_step_decrease_dca} the exact behavior of one iteration of \eqref{eq:DCA_it} in the setting where $f_2$ may be weakly convex, providing a lower bound for the decrease of the objective expressed in terms of differences of subgradients. 

\item \cref{thm:one_step_decrease_dca} describes a total of six distinct regimes, partitioning the parameters space based on smoothness and strong convexity of both functions. We conjecture that these bounds on the objective decrease are tight for all of those six regimes. Among them, only two were previously known and proved by Abbaszadehpeivasti et al. \citeyearpar{abbaszadehpeivasti2021_DCA}, corresponding to the standard \eqref{eq:DCA_it} setting ($f_1$ and $f_2$ convex) where in addition $F$ is required to be both nonconvex and nonconcave.

\item \cref{thm:dca_rates_N_steps} proves that, in our setting allowing $f_2$ weakly convex, DCA converges sublinearly to critical points, with a $\mathcal{O}(\frac{1}{N})$ rate after $N$ iterations, again with six distinct regimes. 
Based on strong numerical evidence, we conjecture that three of those rates are exact for any number of iterations.

\item We show that a split of the objective $F$ allowing weak convexity of $f_2$ can yield better rates than the standard DCA. Moreover, when both functions are smooth, a well-chosen DC splitting may surpass the celebrated gradient descent.

\item As a direct consequence of our in-depth analysis, we can readily transfer the rates of specific regimes to the PGD setting. 

\end{itemize}
We provide a \href{https://github.com/teo2605/DCA_AISTATS25}{GitHub repository} to support the numerical conjectures and to reproduce all the simulations.
\section{\MakeUppercase{Theoretical background}}%
\begin{definition}\label{def:upper_lower_curvature} 
Let $L>0$ and $\mu \le L$. We say that a proper, lower semicontinuous function $f:\mathbb{R}^d\rightarrow \mathbb{R}$ belongs to the class $\mathcal{F}_{\mu, L}(\mathbb{R}^d)$ (or simply $\mathcal{F}_{\mu, L}$) if and only if it has both (i) upper curvature $L$, meaning that $ \tfrac{L}{2} \|\cdot\|^2 - f$ is convex, and (ii) lower curvature $\mu$, meaning that $ f - \tfrac{\mu}{2} \|\cdot\|^2$ is convex. We also define the class $\mathcal{F}_{\mu, \infty}(\mathbb{R}^d)$ which requires only lower curvature $\mu$.
\end{definition}%

Intuitively, the curvature bounds $\mu$ and $L$ correspond to the minimum and maximum eigenvalues of the Hessian for a function $f \in \mathcal{C}^2$. Functions in $\mathcal{F}_{\mu, L}$ must be smooth when $L<\infty$, while $\mathcal{F}_{\mu, \infty}$ also contains nonsmooth functions.
Depending on the sign of the lower curvature $\mu$, a function $f \in \mathcal{F}_{\mu,L}$ is categorized as: (i) weakly convex (or hypoconvex) when $\mu<0$, (ii) convex when $\mu=0$ or (iii) strongly convex for $\mu>0$. 

The subdifferential of a proper, lower semicontinuous convex (l.s.c.) function $f:\mathbb{R}^d \rightarrow \mathbb{R}$ at a point $x \in \mathbb{R}^d$ is defined as
$$
    \partial f(x)  {}\coloneqq{}  \{ g \in \mathbb{R}^d \,|\, f(y) {}\geq{} f(x) + \langle g, y-x \rangle  \ \forall y\in \mathbb{R}^d\}.
$$%
For weakly convex functions, the subdifferential can be defined as follows Bauschke et al. \citeyearpar[Proposition 6.3]{Bauschke_generalized_monotone_operators_2021}. Let $f \in \mathcal{F}_{\mu,\infty}$ be a weakly convex function with $\mu < 0$. Then $\tilde{f}(x) {}\coloneqq{} f(x) - \mu \tfrac{\|x\|^2}{2}$ is convex with a well-defined subdifferential $\tilde{f}(x)$, and we let $\partial f(x)  {}\coloneqq{}  \{ g - \mu x \,|\, g \in \partial \tilde{f}(x) \}$. Finally, if $f$ is differentiable at $x$, then $\partial f(x) = \{ \nabla f(x)\}$.%

\begin{assumption}[Objective and parameters]\label{assumption:curvatures_on_F}
The objective function $F$ in \eqref{eq:DC_program} is lower bounded and can be written $F = f_1 - f_2$, where $f_1 \in \mathcal{F}_{\mu_1,L_1}$ and $f_2 \in \mathcal{F}_{\mu_2,L_2}$, with parameters $\mu_1 \in [0, \infty)$, $L_1 \in (0,\infty]$, $\mu_2 \in (-\infty,\infty)$ and $L_2 \in (\mu_2, \infty]$, such that $\mu_1<L_1$ and $\mu_2<L_2$.
\end{assumption}%
\cref{assumption:curvatures_on_F} runs throughout the rest of this paper and it implies $F \in \mathcal{F}_{\mu_1 - L_2 \,,\, L_1 - \mu_2}$. Allowing function $f_2$ to be concave is directly applicable to analyzing the PGD iteration on strongly convex functions with long stepsizes (see \cref{sec:PGD_equivalence}). We also denote $F_{\textit{lo}}  {}\coloneqq{}  \inf_x F$.

The domain and range of function $f:\mathbb{R}^d \rightarrow \mathbb{R}$ are $\dom f {} {}\coloneqq{} {} \{ x \in \mathbb{R}^d: f(x) < \infty \}$ and $\range f {} {}\coloneqq{} {} \{y \in \mathbb{R}: \exists \text{ } x \in \dom f \text{ with } y = f(x)\}$, respectively. The domain and range of the subdifferential are $\dom \partial f {}={} \{ x \in \mathbb{R}^d: \partial f(x) \neq \emptyset\}$ and $\range \partial f {}={} \cup \{\partial f(x): x \in \dom \partial f \}$, respectively. The convex conjugate of a l.s.c. function $f$ is defined as $f^*(y) {} {}\coloneqq{} {} \sup_{x \in \dom f} \{ \langle y , x \rangle - f(x) \}$, where $f^*$ is closed and convex. %

\noindent \textbf{DCA iteration.} 
\begin{align}\label{eq:DCA_it}\tag{$\text{DCA}^{}$}
\begin{aligned}
    \text{1. } & \text{Select } g_2 \in \partial f_2(x). \\
    \text{2. } & \text{Select } x^{+} \in \argmin_{w\in \mathbb{R}^d} \{f_1(w) - \langle g_2, w \rangle\}.
\end{aligned}    
\end{align}%
With an abuse of notation, a more compact definition of the \eqref{eq:DCA_it} iteration is $x^+ \in \partial f_1^* ( \partial f_2 (x) )$.

The optimality condition in the definition of $x^+$ implies the existence of $g_1^+ \in \partial f_1(x^+)$ such that $g_1^+ = g_2$, where $g_2 \in \partial f_2(x)$. This is the only characterization of $x^+$ used in our derivations. Note that the sequence of iterates produced by \eqref{eq:DCA_it} is not unique.

\begin{assumption}\label{assumption:well_definiteness_DCA}
The subdifferentials of $f_1$ and $f_2$ satisfy the following conditions: $\varnothing \neq \dom \partial f_1 \subseteq \dom \partial f_2$ and $\range \partial f_2\subseteq \range \partial f_1$. 
\end{assumption}%
\begin{proposition}
    Under \cref{assumption:well_definiteness_DCA}, the \eqref{eq:DCA_it} iterations are well-defined, meaning there exists a sequence $\{x^k\}$, starting from $x^0 \in \dom \partial f_1$, generated by $x^{k+1} \in \partial f_1^*(\partial f_2(x^k))$. 
\end{proposition}%
%

Tao and An  \citeyearpar{DCA_Trust_region_1998} note that DCA is typically well-defined, as for any l.s.c. function $f$, it holds $\ri (\dom f) \subseteq \dom \partial f \subseteq \dom f$, where $\ri (\dom f)$ is the relative interior of $\dom f$. The potential weak convexity of $f_2$ represents only a curvature adjustment in the subdifferential definition.%

A critical point $x^*$ satisfies $\partial f_2(x^*) \cap \partial f_1(x^*) \neq \emptyset$. When both functions $f_1$ and $f_2$ are smooth, any critical point $x^*$ is clearly stationary as we have $\nabla F(x^*) = \nabla f_1(x^*) - \nabla f_2(x^*) = 0$. If only $f_2$ is smooth, we have $\partial F(x^*) = \partial f_1(x^*) - \nabla f_2(x^*)$ (Rockafellar and Wets \citeyearpar[Exercise 10.10]{RockWets98}) and criticality also implies stationarity, since $0 \in \partial F(x^*)$. However, if only $f_1$ is smooth, we can only guarantee the inclusion $\partial (-f_2) (x^*) \subseteq -\partial f_2(x^*)$ (Rockafellar and Wets \citeyearpar[Corollary 9.21]{RockWets98}), implying only $\partial F(x^*) \subseteq \nabla f_1(x^*) - \partial f_2(x^*)$, and critical points may not be stationary. 
\begin{proposition}[Sufficient condition for decrease] \label{prop:sufficient_decrease}
Let $f_1 \in \mathcal{F}_{\mu_1,L_1}$ and $f_2 \in \mathcal{F}_{\mu_2,L_2}$. If $\mu_1 + \mu_2 \ge 0$, then the objective function $F = f_1 - f_2$ decreases after each iteration of \eqref{eq:DCA_it}. Moreover, if $\mu_1 + \mu_2 > 0$ that objective decrease is strict, unless $x^+ = x$.
\end{proposition}%
\cref{prop:sufficient_decrease} is proved in \cref{app:proof_prop_2}, inspired by Dinh and Thi \citeyearpar[Theorem 3, Proposition 2]{Dinh_Thi_dca_1997}.%
%
%
\begin{remark}
    Throughout this paper, we assume that the oracle of $\partial f_1^*$ is exact easy to compute, thus we only focus on the progress of the iterations.
\end{remark}
\textbf{Notation:} Superscripts  indicate the iteration index (e.g., $x^k$ represents the $k$-th iterate). 
%
%
%
%
%
%
\section{\MakeUppercase{Convergence analysis}}\label{sec:rates_smooth_case}%
%

%
\begin{theorem}[One-step decrease]\label{thm:one_step_decrease_dca}
Let $f_1 \in \mathcal{F}_{\mu_1,L_1}$ and $f_2 \in \mathcal{F}_{\mu_2,L_2}$ satisfy \cref{assumption:curvatures_on_F,assumption:well_definiteness_DCA}, with at least $f_1$ or $f_2$ smooth, and assume $\mu_1 + \mu_2 > 0$ or $\mu_1=\mu_2=0$. Then after one step of \eqref{eq:DCA_it} we have
\begin{align}\label{eq:tight_dist_1_with_sigmas}
\hspace{-.5em}
    F(x)-F(x^+)
    {}\geq{} \sigma_i \tfrac{1}{2} \|g_1 - g_2\|^2 + \sigma_i^+ \tfrac{1}{2} \|g_1^+ - g_2^+\|^2
\end{align}%
with $g_1 \in \partial f_1(x)$, $g_1^+ \in \partial f_1(x^+)$, $g_2 \in \partial f_2(x)$, $g_2^+ \in \partial f_2(x^+)$, and the expressions for $\sigma_i, \sigma_i^+ \geq 0$ correspond to one of the six regimes (indexed by $i = 1,\dots,6$) described in \cref{tab:DCA_regimes_one_step} according to the values of parameters $L_1$, $L_2$, $\mu_1$, $\mu_2$. %
\end{theorem}%
\begin{table*}[t]
\centering
\caption{Exact decrease after one iteration: $F(x)-F(x^+)\geq \sigma_i \frac{1}{2} \|g_1-g_2\|^2 + \sigma_i^+ \frac{1}{2} \|g_1^+ - g_2^+ \|^2$ 
(see \cref{thm:one_step_decrease_dca}), with $\sigma_i, \sigma_i^+ \geq 0$ and $p_i = \sigma_i + \sigma_i^+$, $i$=1,\dots,6. The domains satisfy the condition $\mu_1 + \mu_2 > 0$ or $\mu_1 = \mu_2 = 0$, with at least one between $f_1$ and $f_2$ smooth. Notation: $B  {}\coloneqq{}  \mu_1^{-1} + \mu_2^{-1} + L_2^{-1}$ and $E  {}\coloneqq{}  \frac{L_2+\mu_2}{L_1 L_2} \frac{L_2 - L_1}{-\mu_2} + \mu_1^{-1} - L_1^{-1}$.\\}
\label{tab:DCA_regimes_one_step}
\def\arraystretch{1.3}%
\resizebox{\linewidth}{!}{%
\begin{tabular}{|@{}c@{}|@{}c@{}|@{}c@{}|c|c|c|c|}
\hline
\textbf{Regime} &
  \centered{\textbf{$\sigma_i$}} &
  \centered{\textbf{$\sigma_i^+$}}
& \multicolumn{3}{c|}{\textbf{Domain}} & \textbf{Description} \\ [3pt] \hline
\multirow{2}{*}{\textbf{$p_1$}} & 
\multirow{2}{*}{$L_2^{-1} \dfrac{L_2-\mu_1}{L_1-\mu_1}$} &
\multirow{2}{*}{$L_2^{-1} \bigg( 1 + \dfrac{L_2^{-1}-L_1^{-1}}{\mu_1^{-1} - L_1^{-1} } \bigg)$} &
\multirow{2}{*}{\makecell[c]{$L_1 \geq L_2 \geq \mu_1 \geq 0$; \\ $L_1 > \mu_2$}}    & 
\multicolumn{2}{c|}{$\mu_2 \geq 0$} &     \makecell[c]{ $f_1$, $f_2$ convex \\ $F$ nonconvex-nonconcave}        \\  [8pt] \cline{5-7} 
                  &                      &   & &
                \multicolumn{2}{c|}{\makecell[c]{ $\mu_2 < 0$ and $E{}\leq{}0$}} &      \makecell[c]{$f_1$ strongly convex, $f_2$ nonconvex \\ $F$ nonconvex-nonconcave}     
                \\ [8pt] \hline
\textbf{$p_2$} &
  $L_1^{-1} \bigg( 1 + \dfrac{L_1^{-1}-L_2^{-1}}{\mu_2^{-1} - L_2^{-1}} \bigg) $ &
  $L_1^{-1} \dfrac{L_1-\mu_2}{L_2-\mu_2}$
  &
\makecell{$L_2 \geq L_1 \geq \mu_2 \geq 0$;\\$L_2 > \mu_1$}
&
  \multicolumn{2}{c|}{\makecell[c]{$\mu_1 \geq 0$}}
  & \makecell[c]{$f_1$, $f_2$ convex \\ $F$ nonconvex-nonconcave}  
  \\ [8pt] \hline
\multirow{2}{*}{\textbf{$p_3$}} &
\multirow{2}{*}{$\dfrac{L_1^{-1}  \big( \mu_1^{-1} + \mu_2^{-1} + L_2^{-1} \big) }{ \mu_1^{-1} + \mu_2^{-1} + L_2^{-1} - L_1^{-1} }$} &
\multirow{2}{*}{$\dfrac{1}{L_2+\mu_2}$} &
\multirow{2}{*}{\makecell[c]{$\mu_2 < 0$, $\mu_1>0$; \\
$L_2 > \mu_1$; $L_1 > \mu_2$}}
& \multirow{2}{*}{\makecell[c]{$B {}\leq{} 0$}} &
  \makecell{
  $L_1 \geq L_2$ and $E {}\geq{} 0$
  }
  & \multirow{2}{*}{\makecell[c]{$f_1$ strongly convex, $f_2$ nonconvex \\ $F$ nonconvex-nonconcave}}  
  \\ [8pt] \cline{6-6}
& & & & & 
  \makecell{
  $L_2 > L_1$
  }
  & 
  \\ [8pt] \hline
\multirow{3}{*}{\textbf{$p_4$}} &
 \multirow{3}{*}{ $0$}
  &
  \multirow{3}{*}{$\dfrac{\mu_1+\mu_2}{\mu_2^2}$} &
   \multirow{3}{*}{\makecell[c]{$\mu_2 < 0$, $\mu_1 > 0$; \\ $L_1 > \mu_2$; \\
    }}
& \multirow{2}{*}{$B > 0$} & $L_2 > \mu_1 > 0$ &  \makecell[c]{$f_1$ strongly  convex, $f_2$ nonconvex \\ $F$ nonconvex-nonconcave}  
  \\ [8pt] \cline{6-7} 
  & & & & & $0 < L_2 \leq \mu_1$ & \makecell[c]{$f_1$ strongly convex, $f_2$ nonconvex \\ $F$ convex}  
\\ [8pt] \cline{5-7} 
& & & & $B \leq 0$ & $L_2 \leq 0$ & \makecell[c]{$f_1$ strongly convex, $f_2$ concave \\ $F$ strongly convex}  
\\ [8pt] \hline
\multirow{2}{*}{\textbf{$p_5$}} &   
  \multirow{2}{*}{$0$}
  &
  \multirow{2}{*}{$\dfrac{L_2 + \mu_1}{L_2^2}$} &
  \multirow{2}{*}{\makecell[c]{$L_1 > \mu_1 \geq L_2 > 0$;\\$L_1>\mu_2$}}
  & 
  \multicolumn{2}{c|}{$\mu_2 \geq 0$}
  & \makecell[c]{$f_1$ strongly convex, $f_2$ convex \\ $F$ convex}
  \\ [8pt] \cline{5-7} 
  & 
  & & & 
  \multicolumn{2}{c|}{\makecell{$\mu_2 < 0$ and $B \leq 0$}}
  & \makecell[c]{$f_1$ strongly convex, $f_2$ nonconvex \\ $F$ convex}
  \\ [8pt] \hline
\textbf{$p_6$} &
$\dfrac{L_1 + \mu_2}{L_1^2}$
  &
  $0$ &
  \makecell[c]{$L_2 > \mu_2 \geq L_1 > \mu_1$} 
&
  \multicolumn{2}{c|}{\makecell{$\mu_1 \geq 0$}}
  & \makecell[c]{$f_1$ convex, $f_2$ strongly convex \\ $F$ concave}
  \\ [8pt] \hline%
\end{tabular}}
\end{table*}%
\begin{figure*}[t]
    \centering
    \includegraphics[width=\textwidth]{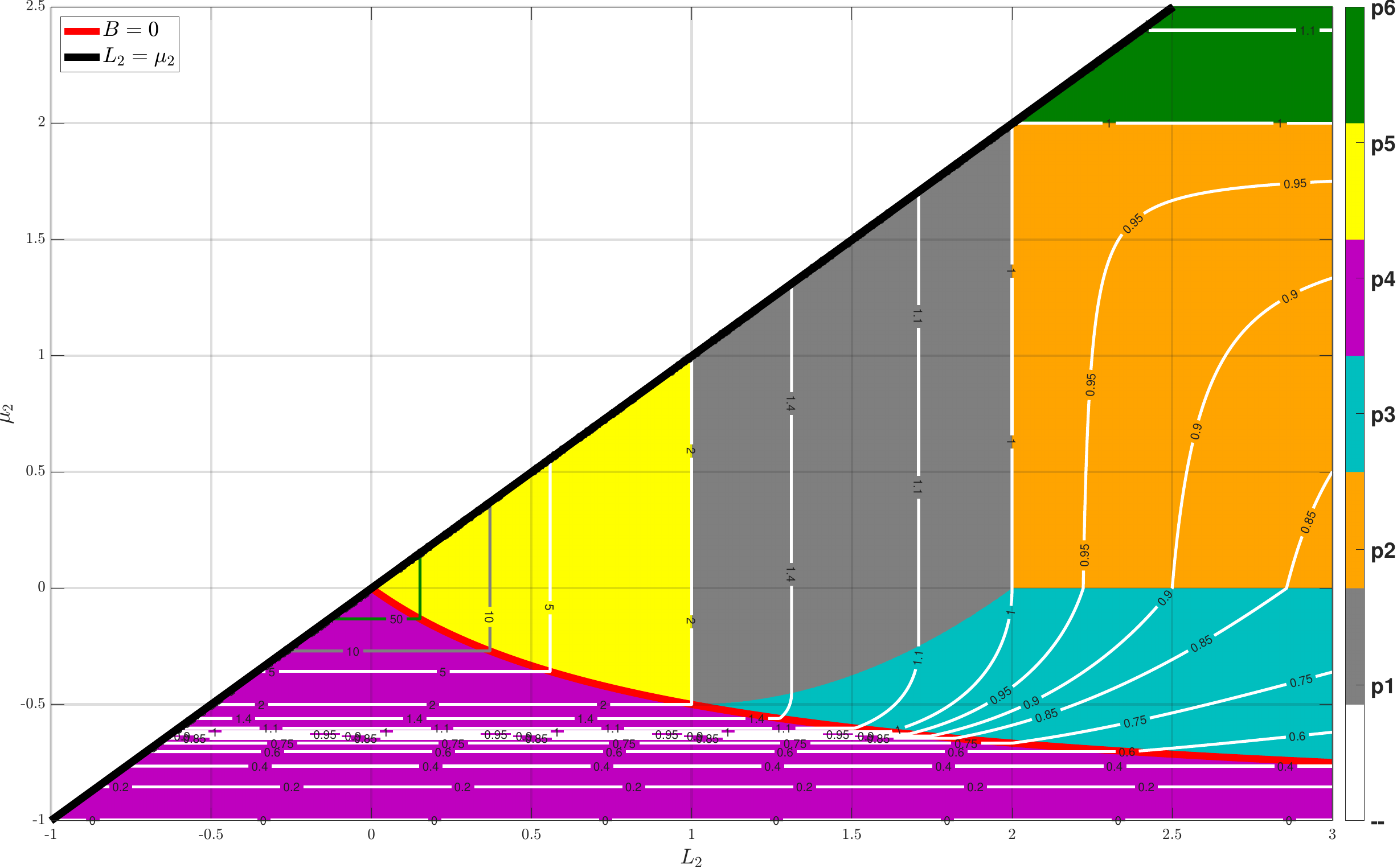}
    \caption{All regimes after one \eqref{eq:DCA_it} iteration (\cref{thm:one_step_decrease_dca}), with $\mu_1 = 1$ and $L_1 = 2$. Contour lines of the denominators $p_i$ vs. $L_2$ and $\mu_2$ are shown. Regimes are bounded by $L_2 > \mu_2$ and $\mu_2 > -\mu_1 = -1$. Regimes $p_1$, $p_2$ and $p_3$ lie within the area delimited by the \textit{threshold} $B=0$ (\eqref{eq:threshold_cond_p3}) and the conditions $L_2 > \mu_1 = 1$ and $L_1 = 2 > \mu_2$ (namely $F$ is nonconvex-nonconcave), and are conjectured to be tight after $N$ iterations of \eqref{eq:DCA_it}.
    }
    \label{fig:all_together}
\end{figure*}%
The six regimes appearing in \cref{tab:DCA_regimes_one_step} are illustrated in \cref{fig:all_together}; we refer to each $p_i$ as one of the six regimes together with its corresponding expression. Notably, there is a striking symmetry between regimes $p_1$ and $p_2$, as well as between $p_5$ and $p_6$. Specifically, the formulas for $p_2$ and $p_6$ in \cref{tab:DCA_regimes_one_step} can be derived from those of $p_1$ and $p_5$ by swapping $L_1 \leftrightarrow L_2$, $\mu_1 \leftrightarrow \mu_2$, and $\sigma_i \leftrightarrow \sigma_i^+$. The proof of \cref{thm:one_step_decrease_dca} is deferred to \cref{sec:proof_one_step_decrease_dca}.%
\begin{conjecture}[Tightest decrease after one iteration]\label{conjecture:tightness_one_N=1}
All six regimes outlined in \cref{thm:one_step_decrease_dca} are tight, i.e., the corresponding lower bounds on the objective decrease cannot be improved.
\end{conjecture}%
\cref{conjecture:tightness_one_N=1} asserts that our set of six inequalities 
represents the tightest possible characterization after a single iteration. Specifically, there exist (separate) function examples for which each of these inequalities holds with equality.%
\begin{corollary}[DCA sublinear rates]\label{thm:dca_rates_N_steps}
Let $f_1 \in \mathcal{F}_{\mu_1,L_1}$ and $f_2 \in \mathcal{F}_{\mu_2,L_2}$ satisfying \cref{assumption:curvatures_on_F,assumption:well_definiteness_DCA}, assume at least $f_1$ or $f_2$ is smooth, and assume $\mu_1 + \mu_2 > 0$ or $\mu_1=\mu_2=0$. Then after $N$ iterations of \eqref{eq:DCA_it} starting from $x^0$ we have
\begin{align}\label{eq:rate_N_steps_no_F*}
    \tfrac{1}{2} \min_{0\leq k \leq N} \{\|g_1^k - g_2^k\|^2\} {}\leq{}
    \frac{F(x^0)-F(x^N)}{p_i(L_1, L_2, \mu_1, \mu_2) N},
\end{align}
where $g_1^k \in \partial f_1(x^k)$ and $g_2^k \in \partial f_2(x^k)$ for all $k=0,\dots,N$ and $p_i=\sigma_i + \sigma_i^+$ is given in \cref{tab:DCA_regimes_one_step}. Additionally, if $F$ is nonconcave (i.e., $L_1 > \mu_2$):
\begin{align*}
    \tfrac{1}{2} \min_{0\leq k \leq N} \{\|g_1^k - g_2^k\|^2\} \leq
    \frac{F(x^0)-F_{\textit{lo}}}{p_i(L_1, L_2, \mu_1, \mu_2) N + \frac{1}{L_1 - \mu_2}}.
\end{align*}
\end{corollary}%
Regimes $p_1$ and $p_2$ correspond in part to the standard setting of \eqref{eq:DCA_it}, where both functions are convex ($\mu_1 \ge 0$, $\mu_2 \ge 0$). Whether $p_1$ or $p_2$ holds depends on which is larger among $L_1$ and $L_2$. These regimes require the objective $F \in \mathcal{F}_{\mu_1-L_2,\ L_1-\mu_2}$ to be nonconvex ($L_2>\mu_1$) and nonconcave ($L_1 > \mu_2$), and were first established by Abbaszadehpeivasti et al. \citeyearpar{abbaszadehpeivasti2021_DCA}, using performance estimation. 
All other described regimes are novel. %
%
\begin{remark}
In the specific convex scenario $\mu_1 = \mu_2 = 0$, both regimes $p_1$ and $p_2$ hold, as outlined by Abbaszadehpeivasti et al. \citeyearpar[Corollary 3.1]{abbaszadehpeivasti2021_DCA} and the one-step decrease is given by:
$F(x)-F(x^+)\geq \frac{1}{2L_1} \|g_1 - g_2\|^2 + \frac{1}{2L_2} \|g_1^+ - g_2^+\|^2$. The same result is obtained using the Bregman proximal point algorithm perspective by  Faust et al.\citeyearpar[Section 4.2]{bregman_DCA_2023}.
\end{remark}%

If $f_1$ is strongly convex, \cref{thm:one_step_decrease_dca} actually extends regime $p_1$ beyond the difference-of-convex case, i.e., to situations where $f_2$ is weakly convex, such that $\mu_1 > 0 > -\mu_2$. This is valid up to a certain threshold determined by the sign of $E  {}\coloneqq{}  \frac{L_2+\mu_2}{L_1 L_2} \frac{L_2 - L_1}{-\mu_2} + \mu_1^{-1} - L_1^{-1}$. For $p_1$, the condition $E < 0$ holds, while for $E \geq 0$ regime $p_3$ emerges. Moreover, for $L_2 \geq L_1$ it always holds $E \geq 0$. Additionally, the boundary of regime $p_3$ is constrained by the threshold $B \leq 0$, where
\begin{align}\label{eq:threshold_cond_p3} 
        B {} {}\coloneqq{} {} \mu_1^{-1} + \mu_2^{-1} + L_2^{-1}.
\end{align}%
Regime $p_4$, emerging for $L_2 \cdot B > 0$, includes two cases: 
(i) when $F$ is nonconvex-nonconcave; and
(ii) when $F$ is strongly convex (even containing $f_2$ concave with $L_2 \leq 0$). The threshold condition $B = 0$ (depicted by the red curve from \cref{fig:all_together}) distinguishes regime $p_4$ from $p_3$ and $p_5$. The later are separated by the condition $L_2 = \mu_1$, delineating the cases $F$ nonconvex (for $p_3$) and $F$ (strongly) convex (for $p_5$), respectively. For completeness of analysis, we also include regime $p_6$, arising for a (strongly) concave objective (and unbounded from below), with $\mu_2 \geq L_1$.

For the particular setup $L_1 = 2$ and $\mu_1 = 1$, in \cref{fig:all_together} we show, as a contour plot, the values of denominators $p_i$ depending on curvatures $L_2$ and $\mu_2$.

Our numerical investigations show that the sublinear rates for $p_{4,5,6}$ in \cref{thm:dca_rates_N_steps} are not tight beyond a single iteration. In the standard case of DCA with $F$ being nonconvex and nonconcave ($\mu_1, \mu_2 < \min \{L_1, L_2\}$), the threshold condition $B = 0$ separates the tight and non-tight regimes in \cref{thm:dca_rates_N_steps}. %
%
\begin{conjecture}[Tightness of sublinear rates]\label{conjecture:tightness_p1_p2_p3_p4}
The DCA rates corresponding to regimes $p_1$, $p_2$ and $p_3$ from \cref{thm:dca_rates_N_steps} are tight for any number of iterations $N$.%
\end{conjecture}%
\cref{conjecture:tightness_p1_p2_p3_p4} asserts that regimes $p_{1,2,3}$ remain tight when exploiting the analysis for one iteration to obtain rates after an arbitrary number of iterations; on such functions, one recovers exactly the worst-case performance when applying DCA. In these cases, closed-from worst-case function examples can be derived. Regime $p_2$ is shown to be tight by Abbaszadehpeivasti et al. \citeyearpar[Example 3.1]{abbaszadehpeivasti2021_DCA} for the specific decomposition $f_1 \in \mathcal{F}_{0,L_1}$ and $f_2 \in \mathcal{F}_{0,\infty}$. In \cref{sec:appendix:tightness}, we provide worst-case examples when both $f_1$ and $f_2$ are smooth, alongside PEP-based numerical evidences.%
\subsection*{One Nonsmooth Term}%
All the above results hold when at least one of the functions $f_1$ and $f_2$ is smooth. When exactly one of them is smooth, i.e., when the other is nonsmooth, some expressions in \cref{tab:DCA_regimes_one_step} become simpler, and we give an explicit description below. In the standard use of DCA, the conjugate step is applied to $f_1$ nonsmooth.%
\begin{corollary}
\label{corollary:dca_rates_nonsmooth_N_steps}
Let $f_1 \in \mathcal{F}_{\mu_1,L_1}$ and $f_2 \in \mathcal{F}_{\mu_2,L_2}$, where exactly one function $f_1$ or $f_2$ is smooth, and assume $\mu_1 + \mu_2 > 0$ or $\mu_1=\mu_2=0$. Consider $N$ iterations of \eqref{eq:DCA_it} starting from $x^0$. Then:%
$$
    \tfrac{1}{2} \min_{0\leq k \leq N} \{\|g_1^k - g_2^k\|^2\} {}\leq{}
    \frac{F(x^0)-F(x^N)}{p_i(L_1, L_2, \mu_1, \mu_2) N},
$$%
where $g_1^k \in \partial f_1(x^k)$ and $g_2^k \in \partial f_2(x^k)$ for all $k=0,\dots,N$ and $p_i$ is provided in \cref{tab:DCA_regimes_nonsmooth}. Additionally, if $F$ is nonconcave (i.e., $L_1 > \mu_2$):%
$$
    \tfrac{1}{2} \min_{0\leq k \leq N} \{\|g_1^k - g_2^k\|^2\} {}\leq{}
    \frac{F(x^0)-F_{\textit{lo}}}{p_i(L_1, L_2, \mu_1, \mu_2) N + \frac{1}{L_1 - \mu_2}}.
$$%
\end{corollary}%
\begin{table}[t]
\centering
\caption{Exact decrease after one iteration: $F(x)-F(x^+)\geq \sigma_i \tfrac{1}{2} \|g_1-g_2\|^2 + \sigma_i^+ \tfrac{1}{2} \|g_1^+ - g_2^+ \|^2$, 
with $f_1$ or $f_2$ nonsmooth, $\sigma_i, \sigma_i^+ \geq 0$, $p_i = \sigma_i + \sigma_i^+$.\\
}
\label{tab:DCA_regimes_nonsmooth}
\resizebox{\linewidth}{!}{%
\begin{tabular}{|c|c|c|c|}
\hline
\textbf{Regime} &
  \textbf{$\sigma_i$} &
  \textbf{$\sigma_i^+$} &
  \textbf{Domain} \\ [2pt] 
  \hline
\textbf{$p_{1,5}$} &
  $0$ &
  $\dfrac{L_2 + \mu_1}{L_2^2} $
  & \makecell{ $L_1 = \infty > L_2 \geq \mu_1 \geq 0$ %
  \\ 
  $\mu_2 \big(\mu_1^{-1} + \mu_2^{-1} + L_2^{-1}\big) \geq 0$
  }
  \\  [6pt] \hline
\textbf{$p_{2,6}$} &
  $\dfrac{L_1 + \mu_2}{L_1^2}$ &
  $0$
  & \makecell{ $L_2 = \infty > L_1 \geq \mu_2 \geq 0$ %
  \\
  $\mu_1 \big(\mu_1^{-1} + \mu_2^{-1} + L_1^{-1}\big) \geq 0$
  }
  \\  [6pt] \hline
\textbf{$p_3$} &
  $\dfrac{\frac{1}{L_1} \big( \frac{1}{\mu_1} + \frac{1}{\mu_2} \big) }{ \frac{1}{\mu_1} + \frac{1}{\mu_2}- \frac{1}{L_1} }$ 
  &
  $0$ 
  & \makecell{ $L_2 = \infty$; $\mu_1 > -\mu_2 > 0$
  }
  \\ [6pt] \hline
\textbf{$p_4$} &
  $0$ 
  &
  $\dfrac{\mu_1+\mu_2}{\mu_2^2}$
  & \makecell{$L_1=\infty$; $\mu_1 > -\mu_2 > 0$ \\
    $0{}<{} \mu_1^{-1} + \mu_2^{-1} + L_2^{-1}$}
  \\  [6pt] \hline
\end{tabular}}%
\end{table}%
%
\cref{corollary:dca_rates_nonsmooth_N_steps} is derived by setting $L_1=\infty$ or $L_2=\infty$ in \cref{thm:dca_rates_N_steps} and in the corresponding entries from \cref{tab:DCA_regimes_one_step}. It shows identical rates as Abbaszadehpeivasti et al. \citeyearpar[Corollary 3.1]{abbaszadehpeivasti2021_DCA} for $\mu_1, \mu_2 \geq 0$, while extending them to scenarios involving weakly convex $f_2$, up to the threshold $B=0$, beyond which regime $p_4$ emerges. Notably, we observe in \cref{tab:DCA_regimes_nonsmooth} that regimes $p_5$ and $p_6$ condense to regimes $p_1$ and $p_2$, respectively. When $L_1=\infty$, only regimes $p_1$ and $p_4$ hold, separated by the threshold $B=0$. Conversely, when $L_2=\infty$, $p_2$ covers the domain with $\mu_2 \geq 0$ and $p_3$ corresponds to $\mu_2 < 0$. For a graphical intuition of these regimes in the nonsmooth case, see \cref{app:sec:nonsmooth_plots}.%
%
\section{\MakeUppercase{Motivations for weakly convex \texorpdfstring{$f_2$}{f2}}}\label{subsec:shifting_curvature}
We assume that the convex conjugate of $f_1$ and of any curvature adjustment $f_1 - \lambda \frac{\|\cdot\|^2}{2}$, where $\lambda \leq \mu_1$, can be computed efficiently. This leads to the natural question: given a splitting $F= f_1 - f_2$, what is the optimal curvature shift $\lambda$ in the decomposition $F = (f_1 - \lambda \frac{\|\cdot\|^2}{2}) - (f_2 - \lambda \frac{\|\cdot\|^2}{2})$? %
In fact, this is a standard approach for addressing weak convexity in the function $f_2$, by lifting it to a convex function, with some $\lambda \leq \mu_2 < 0$, and then applying the DCA iteratively. 

We show that this approach is suboptimal and the sublinear rates, compared in terms of largest denominator $p_i$, can be improved. 
Specifically, the constant in sublinear rates for a nonconvex-nonconcave objective $F$ (where $\max\{\mu_1,\mu_2\} < \min\{L_1,L_2\}$) is optimized with respect to the curvature shift $\lambda$ by maximizing the denominator. In some cases, as indicated in the \textit{Ratio} column, the improvement is significant.

Let $\tilde{f}^{\lambda}_{i} \coloneqq f_{i} - \lambda \frac{\|\cdot\|^2}{2}$, with $i=\{1,2\}$, be the curvature adjusted functions, and $P_{\lambda} \coloneqq P(L_1,\mu_1,L_2,\mu_2,\lambda)$ be the denominator corresponding to one of the six possible regimes $p_i$, determined by the initial splitting curvatures and the parameter $\lambda$. Initially, $\lambda=0$ and the denominator is $P_0$. In \cref{tab:shifting_curvature_examples} we use $p_i$, along with its value, to represent the regime before and after the curvature adjustment. Given the analytical expressions for all six regimes, we can easily numerically compute $\lambda^*=\argmax_{\lambda} P_{\lambda}$.%
\begin{table*}[!ht]
\centering
\caption{Improvement of DCA convergence rates for nonconvex-nonconcave objectives $F$ ($\max\{\mu_1,\mu_2\} < \min\{L_1,L_2\}$) using curvature shifting with optimal $\lambda^*$. Relative improvement (\textit{Ratio}) is defined as $\frac{P(\lambda^*)-P(0)}{P(0)}$, where $P(0)$ is the initial denominator before splitting and $P(\lambda^*)$ is the optimal (largest) one.\\}
\label{tab:shifting_curvature_examples}
\resizebox{.9\linewidth}{!}{%
\begin{tabular}{|c|ccccc|ccc|}
\hline
\textbf{Setup}         & $\mu_1$ & $L_1$  & $\mu_2$ & $L_2$ & $P(0)$ & $\lambda^*$ & \begin{tabular}[c]{@{}l@{}}$P(\lambda^*)$\end{tabular} & \textbf{Ratio} \\ \hline
\multirow{3}{*}{$\mu_1 > \mu_2$} & $0.2$   & $3$    & $0.1$   & $4$   & $p_2 = 0.4535$          &  $0.1221$           & $p_3=0.6031$              & $33 \%$      \\
                                      & $0.2$   & $1000$ & $0.1$   & $3$   & $p_1=0.3255$          &     $0.1494$        & $p_1=p_3=0.358$   &  $10 \%$         \\
& $1$ & $2$ & $0.5$ & $1.5$ & $p_1=0.76$ & $0.6733$ & $p_3=2.0321$ & $167 \%$
                                      \\ \hline
$\mu_1=\mu_2$                    & $1$     & $4$    & $1$     & $3$   & $p_1=0.4583$          & $\mu_1=1$   & $p_1=p_2=0.8333$ & $81 \%$                                                               \\ \hline
\multirow{2}{*}{$\mu_1 < \mu_2$} & $0.1$   & $3$    & $0.2$   & $4$   & $p_2=0.4539$          & $\mu_1=0.1$ & $p_2=0.602$ & $32 \%$         \\
                                      & $0.001$   & $4$  & $0.002$   & $3$ & $p_1=0.4531$          & $\mu_1=0.001$ & $p_1=0.5835$ & $28 \%$               \\ \hline 
\multirow{2}{*}{\makecell[c]{$\mu_1 > 0$; $\mu_1 + \mu_2 > 0$}} & $2$ & $4$ & $-1.75$ & $3$ & $p_4 = 0.4688$ & $-0.4855$ & $p_3=0.52$ & $11\%$ \\ 
& $2.99$ & $4$ & $-2.9$ & $3$ & $p_4 = 0.4994$ & $-0.935$ & $p_4 = 0.5076$ & $1.6\%$ \\
\hline 
\makecell[c]{$\mu_1 > 0$; $\mu_1 + \mu_2 < 0$} & $1$ & $2$ & $-1.5$ & $1.5$ & -- & $-0.6526$ & $p_3=0.8516$ & -- \\ 
\hline 
\end{tabular}}
\end{table*}%
%
%
%
%

The examples in \cref{tab:shifting_curvature_examples} suggest the following observations. First, consider both functions to be (strongly) convex. If $\mu_1 \leq \mu_2$, then the best splitting is achieved by making $\tilde{f}_1$ convex, hence $\lambda = \mu_1$. Notably, when $\mu_1 = \mu_2$, both functions become convex, which surprisingly implies that the initial strong convexity may actually slow down the algorithm. If $\mu_1 > \mu_2$, the optimal splitting occurs when $f_2$ is shifted to a weakly convex function. Furthermore, even when starting with $f_2$ as weakly convex, the best curvature maintains this weak convexity. Additionally, in case of a \textit{bad} decomposition where $\mu_1 + \mu_2 < 0$ (providing no convergence guarantee for DCA iterations), an appropriate  $\lambda$ can ensure convergence. In fact, the optimal denominator is reached for some $\mu_2 < 0$. 

This collection of seemingly surprising results involving weak convexity of $f_2$ is explained by the equivalence with the proximal gradient descent (PGD) (see \cref{sec:PGD_equivalence}) and its faster convergence when using a stepsize larger than the inverse Lipschitz constant. %
%
\paragraph{Benefits on Smooth Functions.} 
The objective $F$ is smooth when both $f_1$ and $f_2$ are smooth, raising the question of why to not use gradient descent (GD) directly. 
We compare the rates of DCA, having the iteration $x^+ = \nabla f_1^*(\nabla f_2(x))$, to GD, whose iteration with stepsize $\gamma \in (0, \frac{2}{L_F})$ reads $x^+ = x - \gamma \nabla F(x)$. We use the criteria $\tfrac{1}{2}\min\{\|\nabla F(x)\|^2, \|\nabla F(x^+)\|^2\} {}\leq{} \frac{F(x)-F(x^+)}{p}$, where $p$ denotes the worst-case denominator for DCA or GD. In \cref{example:F_smooth_nonconvex} we show that, for smooth nonconvex functions, the optimal rate for DCA applied can surpass the optimal rate of GD.%
\begin{example}[$F$ smooth-nonconvex]\label{example:F_smooth_nonconvex}
Let $L_F = -\mu_F= 1$, and the decomposition $F=f_1-f_2$, with parameters $\mu_1 = 1.5$, $L_1 = 2$, $\mu_2 = 1$ and $L_2 = 2.5$. This corresponds to regime $p_2$, with $p_2=0.9167$. The optimal stepsize for GD, $\gamma^*=\frac{2}{\sqrt{3}}$ (Abbaszadehpeivasti et al. \citeyearpar[Theorem 3]{abbaszadehpeivasti2021GM_smooth}), yields the denominator $p_{\text{GD}}=1.5396$, which is $67\%$ larger than for the initial DCA splitting. However, the best DCA splitting achieved by subtracting curvature $\lambda^*=1.0091$, leading to $F=\tilde{f}_1^{\lambda} -\tilde{f}_2^{\lambda}$ with $\tilde{f}_1^{\lambda} \in \mathcal{F}_{0.4009, 0.9009}$ and $\tilde{f}_2^{\lambda} \in \mathcal{F}_{-0.0991,1.4009}$, corresponds to regime $p_3$, with $p_{\text{DCA*}}=1.724$, which is $12\%$ larger than $p_{\text{GD}}$. Thus, with appropriately chosen curvatures in DCA, we can improve convergence rates even in the smooth case.
\end{example}%
%
%
\section{\MakeUppercase{Improved convergence rates of PGD}}\label{sec:PGD_equivalence}%
In this section we demonstrate that convergence rates for PGD can be directly derived from its iterate equivalence with DCA, a connection that is often underemphasized in the literature. While DCA is typically applied to nonconvex-nonconcave objective functions, it can also serve as a valuable tool to establish rates in the cases involving (strongly) convex costs.%
\begin{assumption}[PGD splitting]\label{ass:PGD_splitting}
    Consider the composite objective function $F=\varphi + h$, where $\varphi \in \mathcal{F}_{\mu_{\varphi}, L_{\varphi}}$ is smooth, with $L_{\varphi}>0$ and $\mu_{\varphi} \in (-\infty, L_{\varphi})$, and $h: \mathbb{R}^d \rightarrow (-\infty, \infty]$ is proper, closed and convex, $h \in \mathcal{F}_{\mu_{h}, L_{h}}$, such that $0 \leq \mu_h \leq L_h$.%
\end{assumption}%
%
%
The PGD iteration with stepsize $\gamma > 0$ is given by:
\begin{align}\label{eq:PGD_iteration}\tag{PGD}
\hspace{-.5em}
\begin{aligned}
    & x^+ {} {}\coloneqq{} {} \prox_{\gamma \varphi} \left[ x - \gamma \nabla \varphi(x) \right] \\
        &{}={} \argmin_{w \in \mathbb{R}^d} \big\{ h(w) + \frac{1}{2\gamma} \left\|w - x + \gamma \nabla \varphi(x) \right\|^2 \big\}.
\end{aligned}
\end{align}
To the best of our knowledge, no tight rate expressions exist for the case when $\varphi$ is nonconvex. For $\varphi$ convex, refer to Taylor et al. \citeyearpar{Taylor_Jota_PGM_rates_proofs}, employing different performance metrics.
\begin{proposition}[PGD is equivalent to DCA; Le Thi and Pham Dinh \citeyearpar{LeThi_2018_30_years_dev}, Section 3.3.4]\label{prop:Equivalence_PGD_DCA}
Starting from $x \in \mathbb{R}^d$, one iteration of \eqref{eq:PGD_iteration} with stepsize $\gamma>0$ on the composite objective function $F=\varphi + h$, under \cref{ass:PGD_splitting}, yields the same point $x^+$ as a DCA iteration applied to the splitting $F=f_1 - f_2$, where $f_1 = h + \frac{1}{2\gamma} \|\cdot\|^2$ and 
$f_2 = \frac{1}{2\gamma}\|\cdot\|^2 - \varphi$. Furthermore, the curvatures are related as follows: $\mu_1 {}={} \gamma^{-1}+\mu_{h} $; $L_1 {}={} \gamma^{-1}+L_{h}$; 
$\mu_2 {}={} \gamma^{-1}-L_{\varphi}$; $L_2  {}={} \gamma^{-1}-\mu_{\varphi}$.
\end{proposition}%
%
\begin{proposition}\label{prop:equiv_performance_metric}
    The convergence measure based on the (sub)gradient residual norm is identical when applying PGD to the decomposition $F = \varphi + h$ or DCA to the decomposition $F = f_1 - f_2$, where $f_1$ and $f_2$ are defined in \cref{prop:Equivalence_PGD_DCA}.
\end{proposition}%
\cref{prop:Equivalence_PGD_DCA,prop:equiv_performance_metric} are proved in \cref{appendix:proofs:Props_equiv_PGD_DCA}. \cref{prop:equiv_performance_metric} shows that for any iteration $x$ it holds $\|\nabla \varphi(x) + g_h\|^2 {}={} \|g_1 - \nabla f_2(x)\|^2$, where $g_h \in \partial h(x)$ and $g_1 \in \partial f_1(x)$.
To determine convergence rates for PGD, one can substitute the curvature values from the DCA convergence rate expressions in \cref{sec:rates_smooth_case} with the corresponding curvatures defined in terms of PGD parameters as specified in \cref{prop:Equivalence_PGD_DCA}. Moreover, the smoothness of $\varphi$ implies the smoothness of $f_2$, ensuring that the critical points are also stationary.

\noindent \textbf{Particular cases of PGD.} 
With $h=0$, PGD reduces to the celebrated gradient descent. If $h = \delta_C$, the indicator function of a non-empty, closed and convex set $C$, then PGD becomes the projected gradient descent. Setting $\gamma = \frac{1}{L}$ yields the iterative shrinkage thresholding algorithm (ISTA). If $\varphi$ is a constant function, PGD becomes the proximal point method. 

We focus on the typical PGD setup, where $\varphi$ is smooth and $h$ is convex and nonsmooth, thus $\mu_h = 0$, $L_h = \infty$. Hence, the DCA-like curvatures are $\mu_1 {}={} \gamma^{-1}$, $L_1 {}={} \infty$, $\mu_2 {}={} \gamma^{-1} - L_{\varphi}$, $L_2 {}={} \gamma^{-1} - \mu_{\varphi}$. However, if additional information about $h$ is available, it can be similarly incorporated to derive improved rates.

From the equivalence of the curvatures, it follows that large stepsizes $\gamma > \frac{1}{L_{\varphi}}$ correspond to negative $\mu_2$, indicating that $f_2$ is weakly convex. Furthermore, the condition $\mu_1 + \mu_2 > 0$ in the DCA setting, ensuring the decrease in the objective after one iteration (see \cref{prop:sufficient_decrease}), translates to the standard upper bound on the stepsize for PGD, which is $\gamma < \frac{2}{L_{\varphi}}$.

In \cref{tab:PGD_regimes_one_step} we summarize the corresponding DCA curvatures for various notable cases in the PGD setting. For the stepsize $\gamma = \frac{1}{L_{\varphi}}$, commonly used in convergence analysis, we have $\mu_2 = 0$. Moreover, the unusual case where $L_2 \leq 0$, meaning that $f_2$ is concave, corresponds to large stepsizes $\gamma \in \big[ \frac{2}{L_{\varphi} + \mu_{\varphi}}, \frac{2}{L_{\varphi}} \big)$ applied to strongly convex objectives.%
%
%
\begin{table*}[!ht]
\centering
\caption{PGD settings in DCA curvatures: $\mu_h = 0$, $L_h = \infty$ imply $\mu_1 = \gamma^{-1} > 0$ and $L_1 = \infty$ (see \cref{prop:Equivalence_PGD_DCA}) \\}
\label{tab:PGD_regimes_one_step}
\def\arraystretch{1.25}
\begin{tabular}{|c|c|c|c|c|}
\hline
\textbf{Convexity of $\varphi$ } &
  \textbf{Stepsize $\gamma$} &
  \textbf{$\mu_2 = \gamma^{-1} - L_{\varphi}$} &
  \textbf{$L_2 = \gamma^{-1} - \mu_{\varphi}$} &
  \textbf{\makecell{Regime}}
  \\ \hline
\multirow{3}{*}{\makecell{\textbf{nonconvex}\\ $\mu_{\varphi} < 0$}} &
$\gamma \in \big( 0 , \frac{1}{L_{\varphi}} \big)$ &
$ \mu_2 > 0 $ &
\multirow{3}{*}{ $L_2 > \mu_1$ } & 
\multirow{2}{*}{ $p_1$ }
\\
  & $\gamma = \frac{1}{L_{\varphi}}$ & $\mu_2 = 0$ &  &
\\
  & $\gamma \in \big( \frac{1}{L_{\varphi}}, \frac{2}{L_{\varphi}} \big)$ & $\mu_2 < 0$ & & $p_1$ or $p_4$
  \\ \hline
\multirow{3}{*}{\makecell{\textbf{convex}\\ $\mu_{\varphi} = 0$}} &
$\gamma \in \big( 0 , \frac{1}{L_{\varphi}} \big)$ &
$ \mu_2 > 0 $ &
\multirow{3}{*}{ $L_2 = \mu_1$ } &
\multirow{2}{*}{ $p_1 = p_5$ }
\\
  & $\gamma = \frac{1}{L_{\varphi}}$ & $\mu_2 = 0$ & &
\\
  & $\gamma \in \big( \frac{1}{L_{\varphi}}, \frac{2}{L_{\varphi}} \big)$ & $\mu_2 < 0$ & & $p_4$
  \\ \hline
\multirow{4}{*}{\makecell{\textbf{strongly convex}\\ $\mu_{\varphi} > 0$}} &
$\gamma \in \big( 0 , \frac{1}{L_{\varphi}} \big)$ &
$ \mu_2 > 0 $ &
\multirow{3}{*}{ $0 < L_2 < \mu_1$ } &
\multirow{2}{*}{ $p_5$ }
\\
  & $\gamma = \frac{1}{L_{\varphi}}$ & $\mu_2 = 0$ & & 
\\
  & $\gamma \in \big( \frac{1}{L_{\varphi}}, \frac{2}{L_{\varphi} + \mu_{\varphi}} \big)$ & \multirow{2}{*}{$\mu_2 < 0$} & &
  $p_4$ or $p_5$
  \\ 
  & $\gamma \in \big[ \frac{2}{L_{\varphi} + \mu_{\varphi}}, \frac{2}{L_{\varphi}} \big)$ & & $L_2 \leq 0 < \mu_1$ & $p_4$
  \\ \hline
\end{tabular}
\end{table*}
%
\begin{theorem}[One-step decrease PGD]\label{thm:one_step_decrease_pgd}%
Let $\varphi \in \mathcal{F}_{\mu_{\varphi},L_{\varphi}}$ and $h \in \mathcal{F}_{0,\infty}$, with $\varphi$ smooth. Consider one iteration of \eqref{eq:PGD_iteration} with stepsize $\gamma \in \big(0, \frac{2}{L_{\varphi}}\big)$. Then:%
\begin{align*}
    F(x)-F(x^+)
    {}\geq{} \sigma^{+} ( L_{\varphi}, \mu_{\varphi}, \gamma ) \tfrac{1}{2} \|\nabla \varphi(x^+)+g_h^+\|^2, 
\end{align*}%
with $g_h^+ \in \partial h(x^+)$ and $\sigma^{+} ( L_{\varphi}, \mu_{\varphi}, \gamma ) \geq 0$ defined as:%
\vspace{-0.5em}
\begin{align*}
    \hspace{-.5em}
    \left\{
    \def\arraystretch{1.5} \arraycolsep=1pt
    \begin{array}{cl}
       \frac{\gamma ( 2 - \gamma \mu_{\varphi} )}{(1 - \gamma \mu_{\varphi})^2},  &
       \text { if }
       \mu_{\varphi} \leq 0 
       \text { and either }
       \hspace{-0.3em}
       \left\{
       \hspace{-0.2em}
       \begin{array}{l}
          \gamma \in \big( 0 , \frac{1}{L_{\varphi}} \big];
          \\
          \makecell[l]{\gamma \in \big( \frac{1}{L_{\varphi}}, \frac{2}{L_{\varphi}} \big)
          \\ \,\, \text { and } B \leq 0;}
       \end{array}\right.\\
       \frac{ \gamma ( 2-\gamma L_{\varphi} )}{(1-\gamma L_{\varphi})^2},  & 
       \text { if } \,\,
       \makecell[l]{\gamma \in \big( \frac{1}{L_{\varphi}}, \frac{2}{L_{\varphi}} \big) 
       \text{ and } B > 0,}
    \end{array}
    \right.
\end{align*}
where $B = 1 + \frac{1}{1-\gamma L_{\varphi}} + \frac{1}{1-\gamma \mu_{\varphi}}$.
\end{theorem}%
\cref{thm:one_step_decrease_pgd} is derived by substituting the curvature expressions from the DCA splitting, as outlined in \cref{prop:Equivalence_PGD_DCA}, into \cref{thm:one_step_decrease_dca}.%
%
%
\begin{corollary}[Sublinear PGD Rates]\label{thm:pgd_rates_N_steps}
Let $\varphi \in \mathcal{F}_{\mu_{\varphi},L_{\varphi}}$ and $h \in \mathcal{F}_{0,\infty}$, with $\varphi$ smooth, as in \cref{ass:PGD_splitting}. Consider $N$ iterations of \eqref{eq:PGD_iteration} with stepsize $\gamma \in \big(0, \frac{2}{L_{\varphi}}\big)$, starting from $x^0$. Then:
$$
    \tfrac{1}{2} \min_{0\leq k \leq N} \{\|\nabla \varphi(x^k) + g_h^k \|^2\} {}\leq{}
    \frac{F(x^0)-F_{\textit{lo}}}{\sigma^{+}(L_1, L_2, \mu_1, \mu_2) N},
$$%
where $g_h^k \in \partial h(x^k)$ for all $k=1,\dots,N$ and $\sigma^{+}(L_1, L_2, \mu_1, \mu_2) \geq 0$ is given in \cref{thm:one_step_decrease_pgd}.
\end{corollary}%
\cref{thm:pgd_rates_N_steps} follows immediately from \cref{corollary:dca_rates_nonsmooth_N_steps}. \cref{conjecture:tightness_p1_p2_p3_p4} on tightness implies that only the first branch in \cref{thm:pgd_rates_N_steps} leads to an exact rate. This corresponds to non-strongly convex functions $\varphi$, with $\mu_{\varphi} \leq 0$, and stepsizes $\gamma$ less than some threshold $\bar{\gamma}$ which cancels the expression of the threshold $B$. Further analysis is needed to complete the picture of exact rates for any stepsize $\gamma$ and any curvature of $\varphi$.%
\section{\MakeUppercase{Proof technique}}\label{sec:proofs}%
\textbf{Performance estimation problem (PEP) and interpolation.} 
PEP, introduced by Drori and Teboulle \citeyearpar{drori_performance_2014} and further refined by  Taylor et al. \citeyearpar{taylor_smooth_2017}, provides a framework for analyzing tight convergence rates in various optimization methods. It casts the search for the worst-case function as an optimization problem itself, using specific inequalities which are necessary and sufficient to interpolate functions from the target class, as in \cref{thm:interp_hypo_characterization_min}. The problem is then reformulated as a convex semidefinite program, numerically solvable. Abbaszadehpeivasti et al. \citeyearpar{abbaszadehpeivasti2021_DCA} provide an in-depth guide on setting the PEP for DCA.%
%
\begin{theorem}[$\mathcal{F}_{\mu,L}$-interpolation] \label{thm:interp_hypo_characterization_min}
Given an index set $\mathcal{I}$, let $\mathcal{T} {}\coloneqq{} \{(x^i, g^i, f^i)\}_{i \in \mathcal{I}} \subseteq \mathbb{R}^d \times \mathbb{R}^d \times \mathbb{R}$ be a set of triplets. There exists $f \in \mathcal{F}_{\mu,L}$, with $L>0$ and $\mu \le L$, such that $f(x^i)=f^i$ and $g^i \in \partial f(x^i)$ for all $i \in \mathcal I$, iff the following interpolation inequality holds for every pair of indices $(i,j)$, with $i,j \in \mathcal{I}$:%
\begin{align}\label{eq:Interp_hypoconvex}\tag{$Q_{i,j}$}
\begin{aligned}
    & f^{i}-f^{j}-\langle g^{j}, x^{i}-x^{j} \rangle 
    \geq
    \frac{1}{2L} \|g^{i}-g^{j}\|^{2} {}+{}
    \\ & {}\qquad\quad{} +
    \frac{\mu}{2L(L-\mu)} \|g^{i}-g^{j} - L (x^i-x^j)\|^2.    
\end{aligned}
\end{align}%
\end{theorem}%
\vspace{-1em}
\cref{thm:interp_hypo_characterization_min} is introduced by Taylor et al. \citeyearpar[Theorem 4]{taylor_smooth_2017} for $\mu \geq 0$ and extended to negative $\mu$ by Rotaru et al. \citeyearpar[Theorem 3.1]{rotaru2022tight}. %
%
The following general performance estimation setup for DCA can be easily integrated in one of the specialized software packages PESTO (in Matlab; Taylor et al. \citeyearpar{PESTO}) or PEPit (in Python; Goujaud et al. \citeyearpar{PESTO_python}):%
\begin{align}\label{PEP:DCA}
    \begin{aligned}    
        \maximize \,\, & 
        \frac{ \frac{1}{2} \min \limits_{0 \leq k \leq N} \{\|g_1^k - g_2^k \|^2 \}} {F(x^0)-F(x^N)}
        \\
        \text { subject to }
            \,\, & 
            \{(x^k,g_{1,2}^k,f_{1,2}^k)\}_{k \in \mathcal{I}} \text{ satisfy \eqref{eq:Interp_hypoconvex}} \\
            \,\, & g_1^{k+1} = g_2^{k} \quad k \in\{0, \ldots, N-1\}.
    \end{aligned}
\end{align}%
The decision variables are $x^k$, $g_1^k$, $g_2^k$, $f_1^k$, $f_2^k$, with $k \in \mathcal{I}=\{0,\dots,N\}$. The numerical solutions of \eqref{PEP:DCA} aided our derivations, enabling verification of the rates from \cref{thm:dca_rates_N_steps}, conjecturing their tightness, and selecting for the proofs the necessary interpolation inequalities from the complete set from \cref{thm:interp_hypo_characterization_min}. 

While employing PEP can sometimes make proofs difficult to comprehend, we offer a clearer demonstration of the various regimes in \cref{sec:proof_one_step_decrease_dca}. This involves combining interpolation inequalities for consecutive iterations, with each of the six regimes proof requiring a distinct weighted combination of them.

\cref{conjecture:tightness_one_N=1,conjecture:tightness_p1_p2_p3_p4}
are supported by extensive numerical evidence obtained through the solution of a large number of performance estimation problems \eqref{PEP:DCA}, spanning the range of parameters $L_1, L_2, \mu_1, \mu_2$. For \cref{conjecture:tightness_one_N=1}, addressing tightness after one iteration, \cref{fig:all_together} is generated by fixing $\mu_1$ and $L_1$, creating a sampling grid for $\mu_2$ and $L_2$, and solving the instances of \eqref{PEP:DCA} for $N=1$. The results align with the \textit{exact} values of $p_{1,\dots,6}$ in their respective regions. \cref{conjecture:tightness_p1_p2_p3_p4}, regarding tightness of regimes $p_1$, $p_2$ and $p_3$ after $N$ iterations, was established through a similar sampling grid with $N\leq 20$, confirming that these regimes yield tight sublinear rates.

%
\section{\MakeUppercase{Experiments: Sparse PCA}}\label{sec:Numerical_experiments}
We provide numerical evidence supporting the efficiency of curvature shifting, aligning with the optimal splittings proposed in \cref{subsec:shifting_curvature}, by exploring the sparse principal component analysis (SPCA) problem. Following  Journ\' ee et al. \citeyearpar{JNRS10_Sparse_PCA} and Themelis et al. \citeyearpar{THP20_DC_envelope_Themelis}, we include elastic-net regularization and solve
$$
    \minimize_{x \in \bar{B}(0,1)} F(x)  {}\coloneqq{}  \kappa \|x\|_1 + \eta \frac{\|x\|^2}{2} - \frac{1}{2}x^T \Sigma x,
$$
where $\bar{B}(0,1) {}\coloneqq{} \{x | \|x\|_2\leq 1\}$ is the closed Euclidean unit ball, $\Sigma = A^T A$ the sample covariance matrix, and $\kappa$, $\eta$ are the $l_1$-(sparsity inducing) and $l_2$-regularization parameters, respectively. We denote 
	$f_1(x)  {}\coloneqq{}  \kappa \|x\|_1 + \eta \frac{\|x\|^2}{2} + \delta_{\bar{B}(0,1)}(x)$ and 
	$f_2(x)  {}\coloneqq{}  \frac{1}{2}x^T \Sigma x$, 
with $\delta_C(x)$ as the indicator function of a set $C$. Note that $\mu_1 = \eta$, $L_1=\infty$, $\mu_2 = \min\{\Lambda(\Sigma)\}$ and $L_2 = \max \{\Lambda(\Sigma)\}$, where $\Lambda(\Sigma)$ denotes the eigenvalues of $\Sigma$. Further, we consider the splittings $\tilde{f}^{\lambda}_{i}(x)  {}\coloneqq{}  f_{i}(x) - \lambda \frac{\|x\|^2}{2}$, with $i=\{1,2\}$, parametrized by $\lambda$ (curvature shift); note that $F(x)$ remains unchanged. The subdifferential of the convex conjugate of $\tilde{f}_1^{\lambda}$ can be computed in closed form for any $\lambda \leq \eta$ (see \cref{app:numerical_experiments}). 
We generated a sparse random matrix $A \in \mathbb{R}^{20n \times n}$ ($n=200$), with $10\%$ density, and normalized $\Sigma = AA^T$ by its maximum eigenvalue, obtaining $f_2 \in \mathcal{F}_{\mu_2=0.3882, L_2=1}$. Parameters $\kappa$ and $\eta$ were selected to ensure a sparse but non-trivial solution, and our focus is to demonstrate the impact of curvature shifting on convergence.

We use $M=1000$ random initial points within the unit ball, and only keep the $\bar{M}$ runs converging to the same solution (non-trivial and with desired sparsity). We compare convergence rates for various splittings, namely $\{0, \pm \lambda^*, \pm 0.5 \lambda^* , \lambda_{\text{max}} \}$, where $\lambda^*$ is the optimal curvature shift defined in \cref{subsec:shifting_curvature} and $\lambda_{\text{max}} \coloneq \frac{\mu_1 + \min\{\mu_1, \mu_2\}}{2}$ is the maximum shift guaranteeing convergence. 
We report in \cref{tab:SPCA_experiments_merged_6_eps} the average number of iterations $N_{\varepsilon}$ required to reach a specific accuracy level $\varepsilon \in \{10^{-1}, 10^{-2}, \dots, 10^{-12}\}$ of the squared (sub)gradient norm of $F(x^k)$ (which is equal to $\argmin_{0 \leq k \leq N} \|\tilde{g}_1^k - \nabla \tilde{f}_2^{\lambda}(x^k)\|^2$), for $\varepsilon=\{10^{-1}, 10^{-2}, \dots, 10^{-12}\}$ and $\tilde{g}_1^k \in \partial \tilde{f}_{1}^{\lambda}(x^k)$. Complete tables are provided in \cref{app:numerical_experiments}.%
\begin{table}[!ht]
\centering
\caption{Average number of iterations $N_{\varepsilon}$ for $\kappa=0.02$ and (Case 1) $\eta=0.5 = \mu_1 > \mu_2$; (Case 2) $\eta=0.2=\mu_1 < \mu_2$, for $\{0, \pm \lambda^*, \pm 0.5 \lambda^* , \lambda_{\text{max}} \}$.\\}
\label{tab:SPCA_experiments_merged_6_eps}
\resizebox{\linewidth}{!}{%
    \begin{tabular}{@{}c@{}|@{}r|cccccccc@{}}
\textbf{Case\,} &  \backslashbox[0pt][l]{ $\lambda$}{$\varepsilon$} &  $10^{-2}$ & $10^{-4}$ &   $10^{-6}$ &  $10^{-8}$ &   $10^{-10}$  \\  \hline
\multirow{6}{*}{\textbf{1}} & 0       & 5.74 & 46.70 & 150.17 & 336.27 & 544.63 \\
                        & \,\textcolor{red}{0.4413}  & \textcolor{red}{3.65} & \textcolor{red}{21.49} & \textcolor{red}{70.75}  & \textcolor{red}{154.93} & \textcolor{red}{247.93}  \\
                        & \,-0.2207 & 6.83 & 59.32 & 189.31 & 426.39 & 693.40 \\
                        & \,-0.4413 & 8.05 & 72.30 & 229.13 & 521.24 & 844.69  \\
                        & \,0.2207  & 4.54 & 33.91 & 108.43 & 244.05 & 394.92 \\
                        & \textcolor{blue}{\,0.4441}  & \textcolor{blue}{3.63} & \textcolor{blue}{21.23} & \textcolor{blue}{70.38}  & \textcolor{blue}{153.78} & \textcolor{blue}{246.61} \\ \hline
\multirow{5}{*}{\textbf{2}}                         & 0                          & 5.88 & 49.85 & 153.46 & 339.67 & 547.61  \\
                         & \textcolor{red}{0.2} & \textcolor{red}{4.32} & \textcolor{red}{37.10} & \textcolor{red}{115.90} & \textcolor{red}{256.33} & \textcolor{red}{412.39}  \\
                         & -0.1                       & 6.48 & 55.52 & 172.24 & 382.11 & 616.29 \\
                         & -0.2                       & 7.06 & 61.44 & 192.46 & 426.00 & 686.04 \\
                         & 0.1                        & 5.26 & 43.47 & 135.39 & 298.14 & 479.93 
    \end{tabular}}
\end{table}%

\textbf{Case 1}: $\eta > \mu_2$. Parameters: $\eta =\mu_1 = 0.5$, $\kappa = 0.02$ (400 runs kept out of 1000). Here, $\lambda^* = 0.4413$ and $\lambda_{\text{max}} = 0.4441$. With $\lambda > \mu_2 = 0.3882$, $\tilde{f}^{\lambda}_2$ becomes weakly convex. Using $\lambda^*$ achieves at least a twofold acceleration compared to $\lambda=0$. Higher $\lambda$ values, increasing the nonconvexity of $\tilde{f}^{\lambda}_2$, further improve the convergence. Conversely, adding curvature to both functions ($\lambda<0$) slows the convergence. The improved convergence rates observed when $\tilde{f}^{\lambda}_2$ becomes weakly convex further highlight the significance of studying DCA for weakly convex functions.%

\textbf{Case 2}: $\eta < \mu_2$. Parameters: $\eta=\mu_1=0.2$, $\kappa=0.02$ (703 runs kept out of 1000). Here, $\lambda^* = \lambda_{\text{max}} = \mu_1 = 0.2$. Using $\lambda^*$ to make $\tilde{f}^{\lambda}_1$ convex improves convergence by $\approx 20\%$ compared to the initial splitting. 

In conclusion, $\lambda^*$ from worst-case analysis is also effective for practical curvature shifting, despite better-than-predicted performance of the later. Reducing convexity in both functions improves the rates, similar to larger stepsizes enhancing PGD convergence. 
%
\section{\MakeUppercase{Conclusion}}%
This work thoroughly examines the behavior of a single \eqref{eq:DCA_it} iteration applied to the DC framework extended to accommodate one weakly convex function. We characterized six distinct regimes for the objective decrease based on subgradient differences and conjecture, based on numerical observations, that certain regimes remain tight across multiple iterations. The convergence rate results as a corollary of this analysis. 

The \eqref{eq:DC_program} structure of the objective facilitates curvature shifting by adding/subtracting $\frac{\lambda}{2} \|x\|^2$ to both terms.
However, applying the \eqref{eq:DCA_it} iteration on a split with one weakly convex function can yield better rates than using a modified split that achieves convexity for both functions, as shown by several examples and confirmed by numerical experiments.

We highlight the strong link between PGD and DCA, leveraging their iteration equivalence to translate DCA rates to PGD setups. Further work is necessary for obtaining tight rates for any stepsize and curvature choice in PGD, aided by the simpler DCA analysis.



\subsubsection*{Acknowledgements}
This research was funded within the framework of the Global PhD Partnership KU Leuven - UCLouvain.

\bibliographystyle{apalike} 
\bibliography{ms.bib} 
\section*{Checklist}

 \begin{enumerate}

 \item For all models and algorithms presented, check if you include:
 \begin{enumerate}
   \item A clear description of the mathematical setting, assumptions, algorithm, and/or model. [\textbf{Yes}] 
   \item An analysis of the properties and complexity (time, space, sample size) of any algorithm. [\textbf{Yes}]
   \item (Optional) Anonymized source code, with specification of all dependencies, including external libraries. [\textbf{Yes}; the source code is provided with the supplementary material \href{https://github.com/teo2605/DCA_AISTATS25}{GitHub repository}.]
 \end{enumerate}

 \item For any theoretical claim, check if you include:
 \begin{enumerate}
   \item Statements of the full set of assumptions of all theoretical results. [\textbf{Yes}]
   \item Complete proofs of all theoretical results. [\textbf{Yes}; due to space constraints, they are provided in the appendix.]
   \item Clear explanations of any assumptions. [\textbf{Yes}]
 \end{enumerate}

 \item For all figures and tables that present empirical results, check if you include:
 \begin{enumerate}
   \item The code, data, and instructions needed to reproduce the main experimental results (either in the supplemental material or as a URL). [\textbf{Yes}; the source code to generate our examples, figures and tables is provided with the supplementary material.]
   \item All the training details (e.g., data splits, hyperparameters, how they were chosen). [\textbf{Not Applicable}]
         \item A clear definition of the specific measure or statistics and error bars (e.g., with respect to the random seed after running experiments multiple times). [\textbf{Not Applicable}]
         \item A description of the computing infrastructure used. (e.g., type of GPUs, internal cluster, or cloud provider). [\textbf{Not Applicable}]
 \end{enumerate}

 \item If you are using existing assets (e.g., code, data, models) or curating/releasing new assets, check if you include:
 \begin{enumerate}
   \item Citations of the creator If your work uses existing assets. [\textbf{Not Applicable}]
   \item The license information of the assets, if applicable. [\textbf{Not Applicable}]
   \item New assets either in the supplemental material or as a URL, if applicable. [\textbf{Not Applicable}]
   \item Information about consent from data providers/curators. [\textbf{Not Applicable}]
   \item Discussion of sensible content if applicable, e.g., personally identifiable information or offensive content. [\textbf{Not Applicable}]
 \end{enumerate}

 \item If you used crowdsourcing or conducted research with human subjects, check if you include:
 \begin{enumerate}
   \item The full text of instructions given to participants and screenshots. [\textbf{Not Applicable}]
   \item Descriptions of potential participant risks, with links to Institutional Review Board (IRB) approvals if applicable. [\textbf{Not Applicable}]
   \item The estimated hourly wage paid to participants and the total amount spent on participant compensation. [\textbf{Not Applicable}]
 \end{enumerate}

 \end{enumerate}

\onecolumn

\appendix
\section{PROOFS}

\subsection{Proof of \texorpdfstring{\cref{prop:sufficient_decrease}}{Proposition 2}} \label{app:proof_prop_2}
We use the following Lemma.
\begin{lemma}[Rotaru et al. \citeyearpar{rotaru2022tight}, Lemma 2.5]\label{lemma:hypoconvex_smooth_quad_bounds}
    Let $L > 0$ and $\mu \le L$, and let $f\in \mathcal{F}_{\mu, L}$. Then $\forall x, y \in \mathbb{R}^d$, with $g \in \partial f(y)$, we have the following:
$$
    \tfrac{\mu}{2} \|x-y\|^2 
        {} \leq {} 
    f(x) - f(y) - \langle g, x-y \rangle 
        {} \leq {} 
    \tfrac{L}{2} \|x-y\|^2.
$$
\end{lemma}%

\begin{proof}[Proof of \cref{prop:sufficient_decrease}]  Using \cref{lemma:hypoconvex_smooth_quad_bounds} on function $f_1$ for the pair $(x,y)=(x,x^+)$ and on function $f_2$ with $(x,y)=(x^+,x)$, and summing the inequalities, we obtain:%
$$\begin{aligned}
\frac{\mu_1 + \mu_2}{2} \|x^+ - x\|^2
        {}\leq{} & 
        f_1(x) - f_1(x^+) - \langle g_1^+\,,\, x - x^+ \rangle 
        {}+{} 
        \\ & {}\quad{} 
        f_2(x^+) - f_2(x) - \langle g_2\,,\, x^+ - x \rangle
        {}\leq{}
        \frac{L_1 + L_2}{2} \|x^+ - x\|^2,
\end{aligned}$$%
where $g_1^+ \in \partial f_1(x^+)$ and $g_2 \in \partial f_2(x)$. By using $F(x) = f_1(x) - f_2(x)$ and $g_1^+ = g_2$ we get:%
$$
        \frac{\mu_1 + \mu_2}{2} \|x - x^+\|^2
            \leq
        F(x) - F(x^+) 
        \leq
        \frac{L_1 + L_2}{2} \|x - x^+\|^2,
$$
enough to prove both parts of the proposition.%
\end{proof}%

\subsection{Proof of \texorpdfstring{\cref{thm:one_step_decrease_dca}}{Theorem 2}}\label{sec:proof_one_step_decrease_dca}
In \cref{proofs:tab:DCA_regimes_one_step} we summarize the coefficients for all regimes, along with the corresponding multipliers $\alpha_i$, which are utilized in the subsequent proofs of each regime in part.%
\begin{table}[h]
\centering
\caption{Exact decrease after one step: $F(x)-F(x^+)\geq \sigma_i \frac{1}{2} \|g_1-g_2\|^2 + \sigma_i^+ \frac{1}{2} \|g_1^+ - g_2^+ \|^2$ 
(see \cref{thm:one_step_decrease_dca}), where $\sigma_i, \sigma_i^+ \geq 0$ and $p_i = \sigma_i + \sigma_i^+$, with $i=1,\dots,6$. Scalar $\alpha_i \geq 0$ is a parameter of the proofs.
}
\label{proofs:tab:DCA_regimes_one_step}
\begin{center}
\begin{tabular}{|c|c|c|c|c|}
\hline
\textbf{Regime} &
  \centered{\textbf{$\sigma_i$}} &
  \centered{\textbf{$\sigma_i^+$}} &
  \textbf{$\alpha_i$} 
  \\ [3pt] 
  \hline
\textbf{$p_1$} &
  $L_2^{-1} \dfrac{L_2-\mu_1}{L_1-\mu_1}$ &
  $L_2^{-1} \bigg( 1 + \dfrac{L_2^{-1}-L_1^{-1}}{\mu_1^{-1} - L_1^{-1} } \bigg) $ &
  $\dfrac{\mu_1}{L_2} \dfrac{L_1-L_2}{L_1-\mu_1}$
  \\ [8pt] \hline
\textbf{$p_2$} &
  $L_1^{-1} \bigg( 1 + \dfrac{L_1^{-1}-L_2^{-1}}{\mu_2^{-1} - L_2^{-1}} \bigg) $ &
  $L_1^{-1} \dfrac{L_1-\mu_2}{L_2-\mu_2}$ &
  $\dfrac{\mu_2}{L_1} \dfrac{L_2-L_1}{L_2-\mu_2}$
  \\ [8pt] \hline
\textbf{$p_3$} &
  $\dfrac{L_1^{-1}  \big( \mu_1^{-1} + \mu_2^{-1} + L_2^{-1} \big) }{ \mu_1^{-1} + \mu_2^{-1} + L_2^{-1} - L_1^{-1} }$ &
  $\dfrac{1}{L_2+\mu_2}$ &
  $\dfrac{-\mu_2}{L_2+\mu_2}$ 
  \\ [8pt] \hline
\textbf{$p_4$} &
  $0$ 
  &
  $\dfrac{\mu_1+\mu_2}{\mu_2^2}$
  &
  $\dfrac{\mu_1+\mu_2}{-\mu_2}$
  \\ [8pt] \hline
\textbf{$p_5$} &
  $0$
  &
  $\dfrac{L_2 + \mu_1}{L_2^2}$ &
  $\dfrac{\mu_1}{L_2}$
  \\ [8pt] \hline
\textbf{$p_6$} &
  $\dfrac{L_1 + \mu_2}{L_1^2}$
  &
  $0$ &
  $\dfrac{\mu_2}{L_1}$ 
  \\ [8pt] \hline%
\end{tabular}
\end{center}
\end{table}%

\begin{proof}[Proof of \cref{thm:one_step_decrease_dca}]
    Let $g_1 \in \partial f_1(x)$, $g_1^+ \in \partial f_1(x^+)$,  $g_2 \in \partial f_2(x)$ and $g_2^+ \in \partial f_2(x^+)$, with $g_1^+ = g_2$. We use the notation: $\Delta x  {}\coloneqq{}  x - x^{+}$, $G  {}\coloneqq{}  g_1 - g_2$, $G^+  {}\coloneqq{}  g_1^+ - g_2^+$ and $\Delta F(x)  {}\coloneqq{}  F(x) - F(x^+)$. By writing \eqref{eq:Interp_hypoconvex} for function $f_1$ with the iterates $(x,x^+)$ we obtain: 
    \begin{align}\label{eq:f_1_evals_in_x_x+}
\begin{aligned}
    f_1(x)-f_1(x^{+}) {}-{} &\langle g_1^+, \Delta x \rangle 
    \geq
    \frac{1}{2L_1} \|G\|^{2} {}+{}
    \frac{\mu_1}{2L_1(L_1-\mu_1)} \|G - L_1 \Delta x\|^2
\end{aligned}
\end{align}%
    and for function $f_2$ with the iterates $(x^+,x)$ we get:
    \begin{align}\label{eq:f_2_evals_in_x+_x}
\begin{aligned}
    f_2(x^+)-f_2(x) {}+{} &\langle g_2, \Delta x \rangle 
    \geq
    \frac{1}{2L_2} \|G^+\|^{2} {}+{}
    \frac{\mu_2}{2L_2(L_2-\mu_2)} \|G^+ - L_2 \Delta x\|^2.
\end{aligned}
\end{align}%
    Summing them up and performing simplifications we get:
    \begin{align}\label{eq:consec_funs_it}
    \begin{aligned}
    \Delta F(x)
    {}\geq{}
    \frac{1}{2L_1} \|G\|^{2} + 
    \frac{\mu_1}{2L_1(L_1-\mu_1)} \|G - L_1 \Delta x\|^2 {}+{}
    \frac{1}{2L_2} \|G^+\|^{2} +
    \frac{\mu_2}{2L_2(L_2-\mu_2)} \|G^+ - L_2 \Delta x\|^2.
    \end{aligned}   
    \end{align}
By writing \eqref{eq:Interp_hypoconvex} for function $f_2$ with the iterates $(x,x^+)$:%
\begin{align}
\begin{aligned}
    f_2(x)-f_2(x^+){}-{}&\langle g_2^{+}, \Delta x\rangle 
    \geq
    \frac{1}{2L_2} \|G^{+}\|^{2} {}+{}
    \frac{\mu_2}{2L_2(L_2-\mu_2)} \|G^{+} - L_2 \Delta x\|^2
\end{aligned}
\end{align}%
and summing it up with \eqref{eq:f_2_evals_in_x+_x} we get:%
\begin{align}\label{eq:dist_1_correction}
\hspace{-.9em}
    \langle G^{+}, \Delta x \rangle 
        \geq
    \frac{1}{L_2} \|G^+\|^{2} +
    \frac{\mu_2}{L_2(L_2-\mu_2)} \|G^+ - L_2 \Delta x\|^2.
\end{align}%
Similarly, by writing \eqref{eq:Interp_hypoconvex} for function $f_1$ with the iterates $(x^+,x)$ and summing it up with \eqref{eq:f_1_evals_in_x_x+} we get:%
\begin{align}\label{eq:dist_1_correction_p_even}
\begin{aligned}
     \langle G, \Delta x  \rangle 
    \geq
    \frac{1}{L_1} \|G\|^{2} +
    \frac{\mu_1}{L_1(L_1-\mu_1)} \|G - L_1 \Delta x\|^2.
\end{aligned}
\end{align}%
The proofs only involve adjusting the right-hand side of \eqref{eq:consec_funs_it} using either inequality \eqref{eq:dist_1_correction} or inequality \eqref{eq:dist_1_correction_p_even}, weighted by scalars $\alpha > 0$, which depend on the curvatures (see \cref{proofs:tab:DCA_regimes_one_step}). Specifically, to establish regimes $p_{1,3,4,5}$ we substitute $\alpha$ in:
\begin{align}\label{eq:ineq_param_in_alpha}
    \hspace{-1.75em}
    \begin{aligned}
    \Delta F(x)
    {}\geq{}&
    \frac{1}{2L_1} \|G\|^{2} + 
    \frac{\mu_1}{2L_1(L_1-\mu_1)} \|G - L_1 \Delta x\|^2 {}+{}
    \\
    & {}\, 
    \frac{1+2\alpha}{2L_2} \|G^+\|^{2} {}+{}
    \frac{\mu_2(1+2\alpha)}{2L_2(L_2-\mu_2)} \|G^+ - L_2 \Delta x\|^2
    {}-{} \alpha \langle G^{+}, \Delta x  \rangle
    \end{aligned}   
\end{align}%
and to demonstrate regimes $p_{2,6}$ we plug in $\alpha$ in:%
$$
    \begin{aligned}
    \Delta F(x)
    {}\geq{}&
    \frac{1}{2L_2} \|G^+\|^{2} +
    \frac{\mu_2}{2L_2(L_2-\mu_2)} \|G^+ - L_2 \Delta x\|^2     
    \\ 
    & {}\, \frac{1+2\alpha}{2L_1} \|G\|^{2} {}+{} 
    \frac{\mu_1(1+2\alpha)}{2L_1(L_1-\mu_1)} \|G - L_1 \Delta x\|^2 
    {}-{} \alpha \langle G, \Delta x \rangle.    
    \end{aligned} 
$$%

Since the proof is entirely based on algebraic manipulations, by exploiting the symmetry on the right-hand side of the two inequalities we focus on demonstrating the regimes $p_1, p_3, p_4, p_5$. Regimes $p_2$ and $p_6$ are complementary to $p_1$ and $p_5$, respectively, under the condition $\mu_1 \geq 0$; specifically, $p_1 \leftrightarrow p_2$ and $p_5 \leftrightarrow p_6$ (see also \cref{tab:DCA_regimes_one_step}). Their proofs can be obtained by interchanging: (i) the curvature indices $1$ and $2$; and (ii) $G$ and $G^+$.

\textbf{Regime $p_1$:} It corresponds to $L_1 \geq L_2 > \mu_1 \geq 0$; in particular, $L_1 \geq \max \{L_2,\mu_1,\mu_2\}$. By setting $\alpha = \frac{\mu_1}{L_2} \frac{L_1 - L_2}{L_1 - \mu_1}$ in \eqref{eq:ineq_param_in_alpha}, after simplifications and building the squares we get:%
$$
\begin{aligned}
\Delta F(x)
 \geq
    \frac{L_2-\mu_1}{L_2(L_1-\mu_1)}
    \frac{\|G\|^{2}}{2}
    + \frac{\mu_1}{L_2(L_1-\mu_1)} \frac{\|G - L_2 \Delta x\|^2}{2}
    + \frac{1}{L_2} \bigg[1 + \frac{\mu_1(L_1-L_2)}{L_2(L_1-\mu_1)}\bigg]
     \frac{\|G^+\|^2}{2} {}+{}
    \\ \,\, 
    {}+{} \frac{\mu_1 \frac{L_1}{L_2} \, \mu_2 \big[\mu_1^{-1} + \mu_2^{-1} + L_2^{-1} - L_1^{-1} \big(2+\frac{L_2}{\mu_2}\big) \big]}{(L_1-\mu_1)(L_2-\mu_2)}
    \frac{\|G^+ - L_2 \Delta x\|^2}{2}.
\end{aligned}%
$$%
The weight of $\|G-L_2 \Delta x\|^2$ is positive ($\mu_1 \geq 0$), while the weight of $\|G^+-L_2 \Delta x\|^2$ is positive if either (i) $\mu_2 {}\geq{} 0$ or (ii) $\mu_2 < 0$ and $\mu_1^{-1} + \mu_2^{-1} + L_2^{-1} - L_1^{-1}
            {}\leq{} 
        L_1^{-1} (1+\frac{L_2}{\mu_2})$. Under these two cases, by neglecting both mixed terms we get:%
$$
\begin{aligned}
\Delta F(x)
&\geq
    \frac{L_2-\mu_1}{L_2(L_1-\mu_1)}
    \frac{\|G\|^{2}}{2}
    + \frac{1}{L_2} \bigg(1 + \frac{L_2^{-1}-L_1^{-1}}{\mu_1^{-1}-L_1^{-1}}\bigg) \frac{\|G^+\|^2}{2},
\end{aligned}
$$
with equality only if $G = G^+ = L_2 \Delta x$.%

\textbf{Regime $p_3$:} It corresponds to $L_2\geq \mu_1 \geq 0$ and negativity of the weight of $\|G^+-L_2 \Delta x\|^2$ from the proof of $p_1$, i.e., $\mu_2 < 0$ and $\mu_1^{-1} + \mu_2^{-1} + L_2^{-1} - L_1^{-1}
            {}>{} 
        L_1^{-1} (1+\frac{L_2}{\mu_2})$. Since $L_2 \geq \mu_1 > -\mu_2$, we have $L_2+\mu_2 > 0$. By setting $\alpha = \frac{-\mu_2}{L_2+\mu_2}$ in \eqref{eq:ineq_param_in_alpha}, after simplifications and building up a square including all weight of $L_2 \Delta x$ we get:%
\begin{align*}
\begin{aligned}
\Delta F(x)
 {}\geq{} &
    \frac{L_1^{-1} ( \mu_1^{-1}+\mu_2^{-1}+L_2^{-1} ) }{\mu_1^{-1}+\mu_2^{-1}+L_2^{-1}-L_1^{-1}} \frac{\|G\|^2}{2} +
    \frac{1}{L_2 + \mu_2} \frac{\|G^+\|^2}{2} {}+{}
    \\& {}\qquad{} {}+{} 
    \frac{\mu_1 L_1 \mu_2 ( \mu_1^{-1}+\mu_2^{-1}+L_2^{-1}-L_1^{-1} )}{L_2(L_1-\mu_1)(L_2-\mu_2)}
    \frac{\Big\| \frac{G}{\frac{L_1}{L_2} + \frac{\mu_2(L_1-\mu_1)}{\mu_1(L_2+\mu_2)}} - L_2 \Delta x\Big\|^2}{2}.
\end{aligned}%
\end{align*}
Since $\mu_2 < 0$, the weight of the mixed term is nonnegative only if $\mu_1^{-1}+\mu_2^{-1}+L_2^{-1} - L_1^{-1} \leq 0$. Further on, this implies that the weight of $\|G\|^2$ is nonnegative only if $\mu_1^{-1}+\mu_2^{-1}+L_2^{-1} \leq 0$, which is exactly the threshold condition $B \leq 0$, with $B$ defined in \eqref{eq:threshold_cond_p3}. Finally, by neglecting the mixed term we get:%
$$
\hspace{-.5em}
\begin{aligned}
\Delta F(x)
 {}\geq{}&
    \frac{L_1^{-1} ( \mu_1^{-1}+\mu_2^{-1}+L_2^{-1} ) }{\mu_1^{-1}+\mu_2^{-1}+L_2^{-1}-L_1^{-1}} \frac{\|G\|^2}{2} {}+{}
    \frac{1}{L_2 + \mu_2} \frac{\|G^+\|^2}{2},
\end{aligned}%
$$%
with equality only if $G = L_2 \left[\frac{L_1}{L_2} + \frac{\mu_2(L_1-\mu_1)}{\mu_1(L_2+\mu_2)} \right] \Delta x$.%

\textbf{Regime $p_4$:} With $\mu_2<0$, if $\mu_1^{-1} + \mu_2^{-1} + L_2^{-1} > 0$, then the proof for regime $p_3$ breaks. Then we set $\alpha = \frac{\mu_1+\mu_2}{-\mu_2} > 0$ in \eqref{eq:ineq_param_in_alpha}. After simplifications and building up the squares we get:%
$$
        \begin{aligned}
           \Delta F(x)
 {}\geq{}
    \frac{\mu_1 + \mu_2}{\mu_2^2}
    \frac{\|G^+\|^{2}}{2} {} + {} 
    \frac{1}{L_1-\mu_1}
    \frac{\|  G - \mu_1 \Delta x \|^2}{2}
    {}+{} 
    \frac{\mu_1 L_2}{-\mu_2} \frac{\mu_1^{-1} + \mu_2^{-1} + L_2^{-1}}{L_2-\mu_2}
    \frac{\| G^+ - \mu_2 \Delta x \|^2}{2}.
        \end{aligned}%
$$%
%
%
The second mixed term has a positive weight if $L_2 ( \mu_1^{-1} + \mu_2^{-1} + L_2^{-1} ) > 0$. Note that this condition allows $L_2 \leq 0$, bounded bellow through the necessary condition $L_2 > \mu_2 > -\mu_1$. After disregarding the mixed squares with positive weights, we obtain:%
$$
           \Delta F(x)
 {}\geq{}
    \frac{\mu_1 + \mu_2}{\mu_2^2}
    \frac{\|G^+\|^{2}}{2},
$$%
which holds with equality only if $G = \mu_1 \Delta x$ and $G^+ = \mu_2 \Delta x$.%

\textbf{Regime $p_5$:} Assume $0 < L_2 \leq \mu_1$, i.e., $F$ is (strongly) convex, and $\mu_2 ( \mu_1^{-1} + \mu_2^{-1} + L_2^{-1} ) \geq 0$, where either $\mu_2 \geq 0$ or $\mu_1 > -\mu_2 > 0$. By setting $\alpha = \frac{\mu_1}{L_2}$ in \eqref{eq:ineq_param_in_alpha}, after simplifications and building up the squares we get:%
$$
        \begin{aligned}
    \Delta F(x)
        {}\geq{}
    & \frac{\mu_1 + L_2}{L_2^2}
   \frac{\|G^+\|^{2}}{2} +
   \frac{1}{L_1-\mu_1}
    \frac{\|G - \mu_1 \Delta x\|^2}{2} 
    {}+{}
    \frac{\mu_1 \mu_2}{L_2} \frac{\mu_1^{-1} + \mu_2^{-1} + L_2^{-1}}{L_2-\mu_2}
    \frac{\| G^+ - L_2 \Delta x \|^2}{2},
        \end{aligned}%
$$%
where the mixed term can be disregarded to obtain:%
$$
    \Delta F(x)
        {}\geq{}
\frac{L_2 + \mu_1}{L_2^2}
   \frac{\|G^+\|^{2}}{2},
$$%
which holds with equality only if $G = \mu_1 \Delta x$ and $G^+ = L_2 \Delta x$.%
\end{proof}%

\subsection{Proof of \texorpdfstring{\cref{thm:dca_rates_N_steps}}{Corollary 1} }
\begin{proof}[Proof of \cref{thm:dca_rates_N_steps}]
From \cref{thm:one_step_decrease_dca}, by taking the minimum between the (sub)gradients difference norms in \eqref{eq:tight_dist_1_with_sigmas} we get:%
$$
\begin{aligned}    
    F(x)-F(x^+) & {}\geq{} \sigma_i \, \frac{1}{2} \|g_1-g_2\|^2 + \sigma_i^+ \frac{1}{2} \|g_1^+ - g_2^+ \|^2 
        \\ & 
        {}\geq{} p_i \, \frac{1}{2} \min \{ \|g_1-g_2\|^2,\|g_1^+ - g_2^+ \|^2 \},
\end{aligned}%
$$%
where $p_i = \sigma_i + \sigma_i^+$, for $i=1,\dots,6$, are given in \cref{tab:DCA_regimes_one_step}. The rate \eqref{eq:rate_N_steps_no_F*} results by telescoping the above inequality for $N$ iterations and taking the minimum among all (sub)gradients differences norms. A rate with respect to $F_{\textit{lo}}$ is obtained either by applying the trivial bound $F(x^0)-F(x^N) \geq F(x^0)-F_{\textit{lo}}$ or, if $L_1 > \mu_2$, by using a tighter bound like the one demonstrated by Abbaszadehpeivasti et al. \citeyearpar[Lemma 2.1]{abbaszadehpeivasti2021_DCA}:%
    $$
        F(x^N)-F_{\textit{lo}} {}\geq{} \frac{1}{2(L_1 - \mu_2)} \|g_1^N-g_2^N\|^2.
    $$%
    By incorporating this into the telescoped sum and once again taking the minimum over all (sub)gradients differences norms, we obtain the second rate from \cref{thm:dca_rates_N_steps}. Furthermore, note that a necessary condition for the tightness of these sublinear rates is that the residual gradient norms $\|g_1^k - g_2^k\|$ must be equal for any $k = 0, \dots, N$. %
\end{proof}%

\subsection{Proofs of \texorpdfstring{\cref{prop:Equivalence_PGD_DCA}}{Proposition 3} and \texorpdfstring{\cref{prop:equiv_performance_metric}}{Proposition 4} }\label{appendix:proofs:Props_equiv_PGD_DCA}
\begin{proof}[Proof of \cref{prop:Equivalence_PGD_DCA}] The approach is to reformulate a PGD iteration as a DCA one, starting from the proximal step definition:%
\begin{align*}
    x^{+} {}={}& \prox_{\gamma h} 
        \big[
        x-\gamma \nabla \varphi(x)
        \big] \\ 
        {}={}&
        \argmin_{w \in \mathbb{R}^d} \big\{
        h(w) +
        \frac{1}{2\gamma}{\big\|w - x + \gamma \nabla \varphi (x) \big\|^2} \big\}  \\ 
         {}={}&
        \argmin_{w \in \mathbb{R}^d} \big\{
        h(w) + \frac{\big\|w\big\|^2}{2\gamma}
        + \langle \nabla \varphi(x) - \frac{1}{\gamma} x, w \rangle
        \big\}  \\ 
         {}={}&
        \argmin_{w \in \mathbb{R}^d} \big\{
        \big[h(w) + \frac{\big\|w\big\|^2}{2\gamma}\big]
        - \big \langle \nabla \big[ \frac{\|x\|^2}{2\gamma}  - \varphi(x) \big], w \big \rangle
        \big\}.
\end{align*}
Then, from $f_1 {}={} h + \tfrac{1}{2 \gamma} \|\cdot\|^2$ and $f_2 {}={} \tfrac{1}{2 \gamma} \|\cdot\|^2 - \varphi$, we recover the \eqref{eq:DCA_it} iteration for $f_2$ smooth.
\end{proof}

\begin{proof}[Proof of \cref{prop:equiv_performance_metric}]
    For any $x \in \mathbb{R}^d$, let $g_h \in \partial h(x)$ and define $G \coloneqq \nabla \varphi(x) + g_h \in \partial F(x)$. By \cref{prop:Equivalence_PGD_DCA}, we have that $g_1 \coloneqq g_h + \frac{x}{\gamma} \in \partial f_1(x)$. Moreover, the smoothness of $\varphi$ implies the one of $f_2$. Then the following holds:
    $$
       G {}={}
       \nabla \varphi(x) + g_h 
            {}={}
        \Big[\frac{x}{\gamma} - \nabla f_2(x)\Big]
            {}+{} 
        \Big[ g_1 - \frac{x}{\gamma} \Big]
            {}={}
        g_1 - \nabla f_2(x).
    $$
Consequently, $\|G\|^2 {}={} \|\nabla \varphi(x) + g_h\|^2 {}={} \|g_1 - \nabla f_2(x)\|^2$. Since $F$ is unchanged between the two splittings, the conclusion follows.
\end{proof}

\section{\MakeUppercase{Tightness analysis} (see \texorpdfstring{\cref{conjecture:tightness_one_N=1,conjecture:tightness_p1_p2_p3_p4}}{Conjectures 1 and 2})} \label{sec:appendix:tightness}
\begin{figure}[H]
    \centering
    \includegraphics[width=0.45\linewidth]{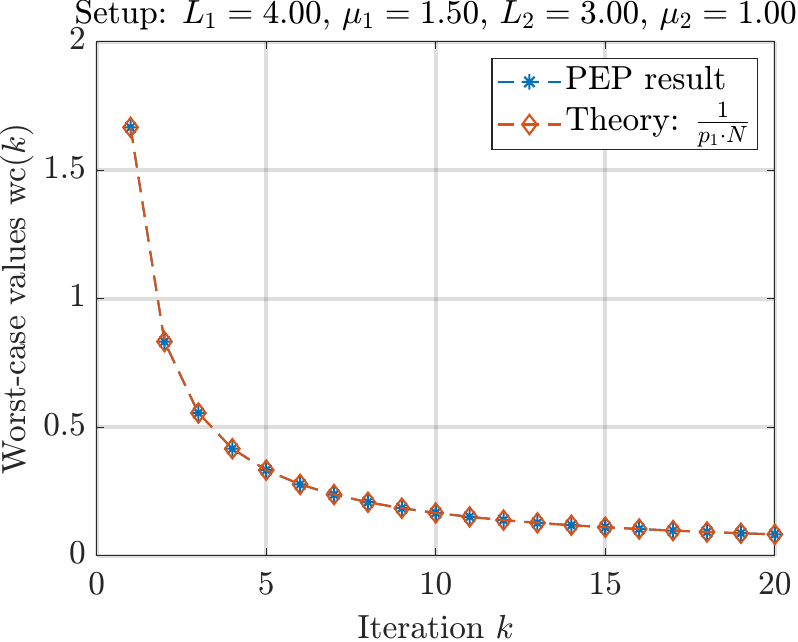}
    \includegraphics[width=0.45\linewidth]{pics/PEP_Sim_Regime_p1_f2_convex.pdf}
    \caption{Examples for regime $p_1$ with $f_2$ convex (\textit{left}) and $f_2$ weakly convex (\textit{right}), showing exactness of our expressions.}
    \label{fig:PEP_regime_p1}
\end{figure}

\begin{figure}[H]
    \centering
    \includegraphics[width=0.45\linewidth]{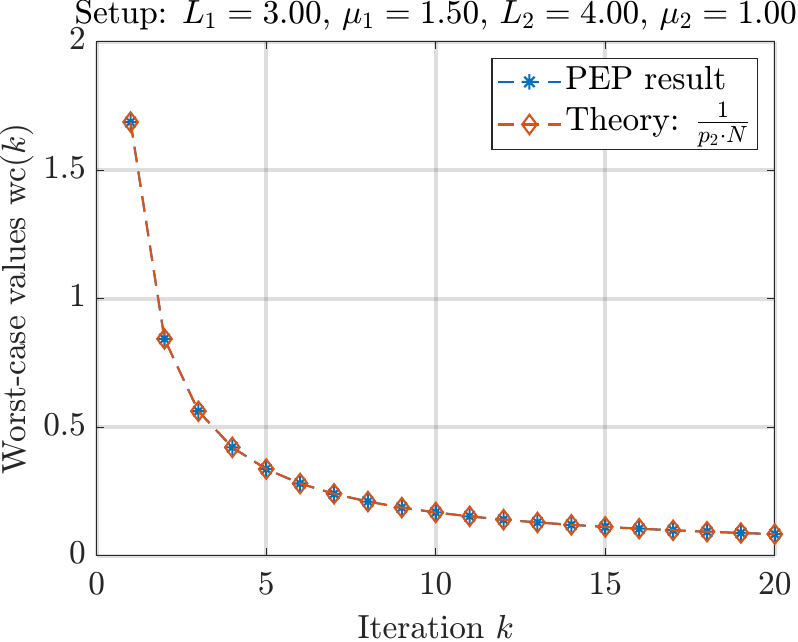}
    \includegraphics[width=0.45\linewidth]{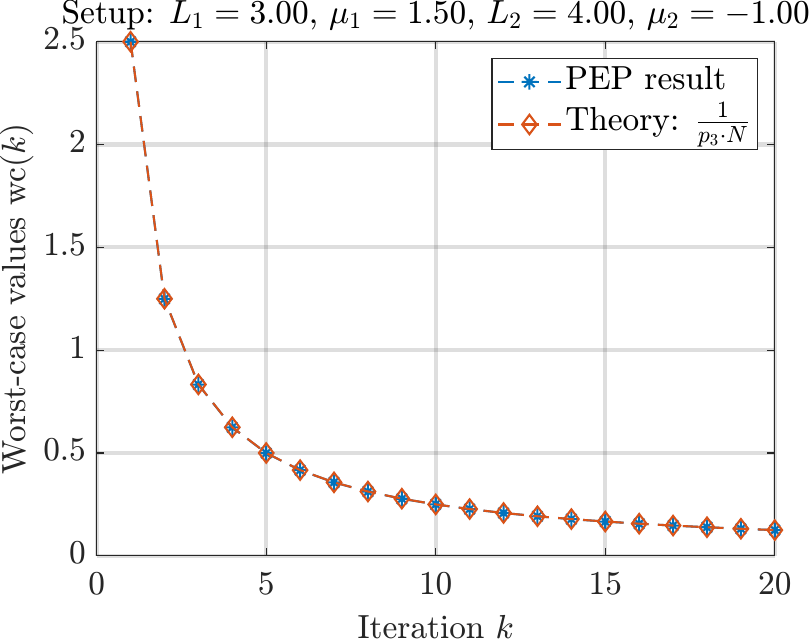}
    \caption{Examples for regimes $p_2$ (\textit{left}) and $p_3$ (\textit{right}), showing exactness of our expressions.}
    \label{fig:PEP_regimes_p2_p3}
\end{figure}

\begin{figure}[H]
    \centering
    \includegraphics[width=0.45\linewidth]{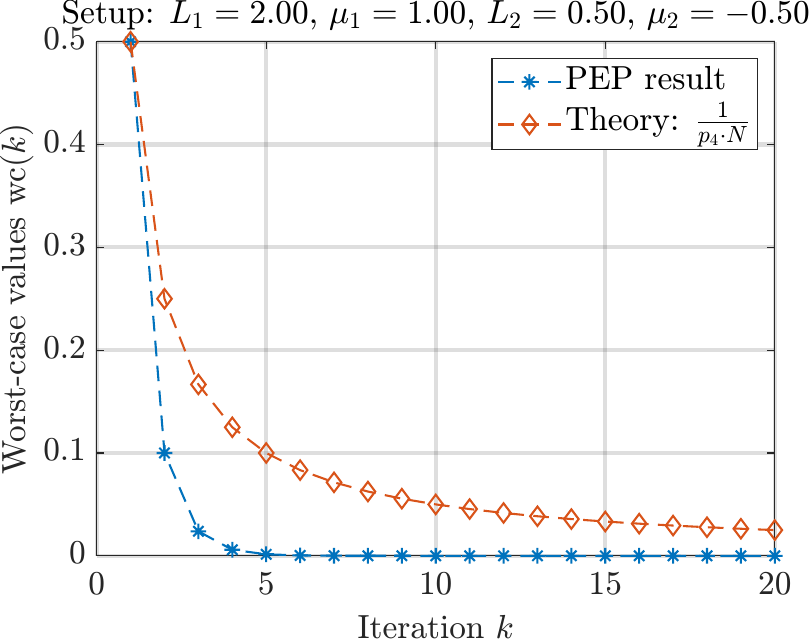}
    \includegraphics[width=0.45\linewidth]{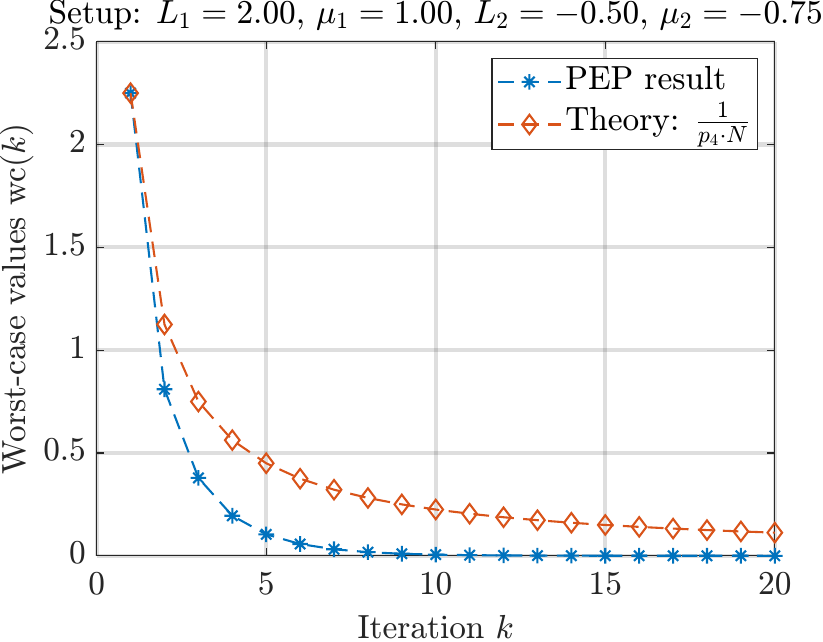}
    \caption{Examples for regime $p_4$, with $L_2>0$ (\textit{left}) and $L_2 < 0$ (\textit{right}), showing that our expressions are non-tight.}
    \label{fig:PEP_regime_p4}
\end{figure}

\begin{figure}[H]
    \centering
    \includegraphics[width=0.45\linewidth]{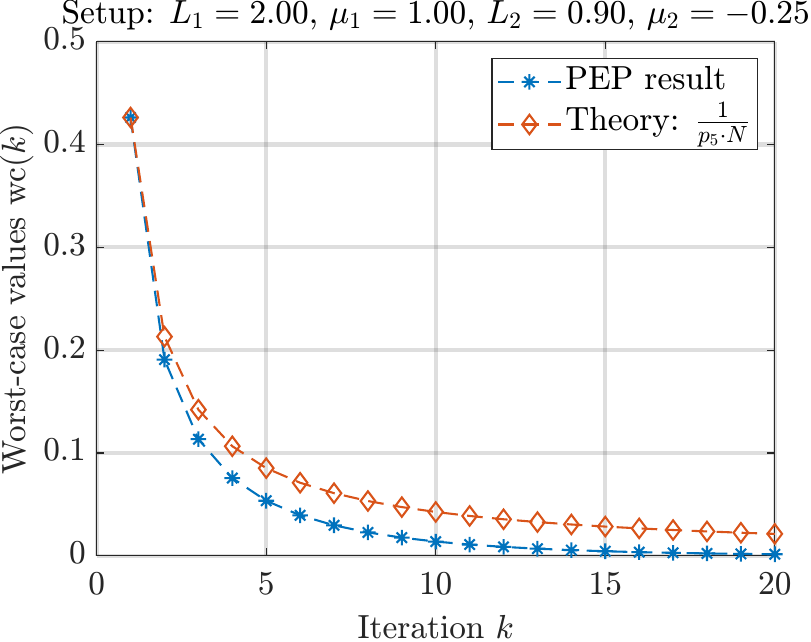}
    \includegraphics[width=0.45\linewidth]{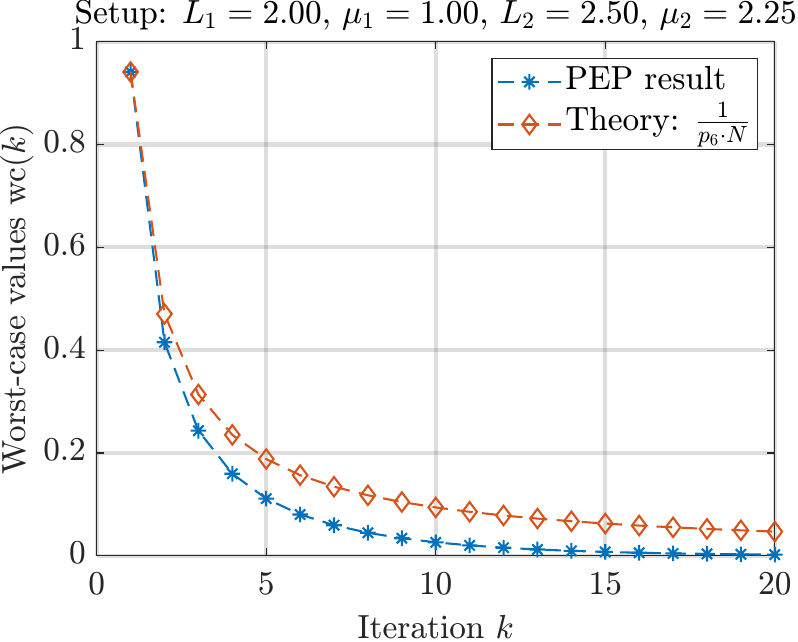}
    \caption{Examples for regimes $p_5$ (\textit{left}) and $p_6$ (\textit{right}), showing that our expressions are non-tight.}
    \label{fig:PEP_regimes_p5_p6}
\end{figure}

For each of the six regimes denoted $p_i$, with $i=1,\dots,6$, we provide numerical examples of the corresponding worst-case values, defined as $\text{wc}(N) = \frac{1}{2} \min_{0 \leq k \leq N} \{\|g_1^k-\nabla f_2(x^k)\|^2\}$, where $g_1^k \in \partial f_1(x^k)$ and the initial condition is fixed to $F(x^0)-F(x^N) \leq 1$. See \cref{fig:PEP_regime_p1} for regime $p_1$ when $f_2$ is convex or weakly convex, \cref{fig:PEP_regimes_p2_p3} for regimes $p_2$ and $p_3$, \cref{fig:PEP_regime_p4} for regime $p_4$ with $L_2>0$ or $L_2<0$ and \cref{fig:PEP_regimes_p5_p6} for regimes $p_5$ and $p_6$. The rates from \cref{thm:dca_rates_N_steps} are predicted to be sublinear, with $\text{wc}(N) = \frac{1}{p_i N}$. This holds for regimes $p_1$, $p_2$ and $p_3$, as stated in \cref{conjecture:tightness_p1_p2_p3_p4}. For regimes $p_4$, $p_5$ and $p_6$, however, the tight rates can be improved and further investigation is required.

The proof of \cref{thm:one_step_decrease_dca} (see \cref{sec:proof_one_step_decrease_dca}) imposes necessary conditions for the worst-cases corresponding to each of the six regimes. Using these, we elaborate on \cref{conjecture:tightness_p1_p2_p3_p4} about the tightness of regimes $p_1$ and $p_2$, providing analytical worst-case functions for the cases when both $f_1$ and $f_2$ are smooth. The nonsmooth case is discussed in Abbaszadehpeivasti et al. \citeyearpar[Example 3.1]{abbaszadehpeivasti2021_DCA} for regime $p_2$, with $f_1 \in \mathcal{F}_{0,L_1}$ and $f_2 \in \mathcal{F}_{0,\infty}$.

\subsection{Worst-Case Function Examples for Regime \texorpdfstring{$p_1$}{p1}}
%
%
\begin{proposition}[Tightness of $p_1$]\label{prop:wc_example_p1}
Let $f_1^0$, $f_2^0$, $g_2^0$, $x^0$ be some real numbers such that $\Delta \coloneqq f_1^0-f_2^0 > 0$ and let $N\geq 1$ be an integer. Consider $\mu_1$, $L_1$, $\mu_2$, $L_2$ as belonging to the domain of regime $p_1$ from \cref{tab:DCA_regimes_one_step}, where the conditions are $0 \leq \mu_1 < L_2 \leq L_1 < \infty$, $\mu_2 < L_1$ and either (i) $\mu_2 \geq 0$ or (ii) $\mu_1 > -\mu_2 > 0$ and $E = \frac{L_2+\mu_2}{L_1 L_2} \frac{L_2 - L_1}{-\mu_2} + \mu_1^{-1} - L_1^{-1} \leq 0$. Let $U{} {}\coloneqq{} {} -\sqrt{\frac{2\Delta}{p_1 N}}$, with expression of $p_1$ given in \cref{tab:DCA_regimes_one_step}. Consider the functions $f_1, f_2: \mathbb{R}\rightarrow \mathbb{R}$ defined as:%
\begin{align*}
    \begin{aligned}
        f_2(x) 
            {} {}\coloneqq{} {}&
        \tfrac{1}{2} L_2 (x-x^0)^2 {}+{}
        g_2^0 (x-x^0) {}+{}
        f_2^0 \\
        f_1(x) {} {}\coloneqq{} {}&
        \left\{
        \def\arraystretch{2}
        \begin{array}{ll}
        \tfrac{1}{2} L_1 (x-x^{0})^2 {}+{}  g_1^{0} (x-x^{0}) {}+{}
            f_1^{0}   & x \in (-\infty, x^0]; \\
            \tfrac{1}{2} L_1 (x-x^k)^2 {}+{}
            g_1^k (x-x^k) {}+{}
            f_1^k  & x \in [x^k, \bar{x}^k];  \\
            \tfrac{1}{2} \mu_1 (x-x^{k+1})^2 {}+{}
            g_1^{k+1} (x-x^{k+1}) {}+{}
            f_1^{k+1}  & x \in [\bar{x}^k, x^{k+1}]; \\
            \tfrac{1}{2} L_1 (x-x^{N})^2 {}+{}
            g_1^{N} (x-x^{N}) {}+{}
            f_1^{N} & x \in [x^N, \infty),
        \end{array}            
        \right.
    \end{aligned}
\end{align*}
with $k=0,1,\dots,N-1$, $\nabla f_2(x^0)=g_2^0$, $f_1(x^k) = f_1^k$, $\nabla f_1(x^k) = g_1^k$ and:%
\begin{align*}
    \begin{aligned}
        x^{k} {}={}& x^0 - k \frac{U}{L_2}, {}\quad{} \forall k = 0,\dots,N; \\
        g_1^{k} {}={}& g_2^0 - (k-1) U, {}\quad{} \forall k = 0,\dots,N; \\
        f_1^{k} {}={}& f_2(x^k) + \frac{N-k}{N} \Delta, {}\quad{} \forall k = 0,\dots,N; \\        
        \bar{x}^{k} {}={}& x^k - \frac{L_2 - \mu_1}{L_1 - \mu_1} \frac{U}{L_2}, {}\quad{} \forall k = 0,\dots,N-1.
    \end{aligned}
\end{align*}
Then by performing $N$ iterations of \eqref{eq:DCA_it} on the function $F(x)=f_1(x)-f_2(x)$, starting from $x^0$, the following result holds:
\begin{align*}
    \tfrac{1}{2} \min_{0\leq k \leq N} \left\{\|\nabla f_1(x^k) - \nabla f_2(x^k)\|^2\right\} {}={}
    \frac{\Delta}{p_1(L_1, L_2, \mu_1, \mu_2) N}.
\end{align*}
\end{proposition}%
\smallskip
\begin{proof}
    By construction, $f_1 \in \mathcal{F}_{\mu_1,L_1}$ and $f_2 \in \mathcal{F}_{\mu_2,L_2}$, where $\mu_2$ can be any real number such that $\mu_2 < L_2$. Moreover, for all $k=0,1,\dots,N-1$ we have $g_2^k {}={} g_1^{k+1}$, which is the necessary condition of solving the DCA iteration. Consequently, the quantity $\|g_1^k - g_2^k\|^2$ is exactly $U^2 = \frac{2 \Delta}{p_1 N}$, where $\Delta = F(x^0) - F(x^N)$. Given that $f_2$ is quadratic and $f_1$ is piecewise quadratic, starting from the initial point $x^0$, the iterations $x^k$ are uniquely determined by applying \eqref{eq:DCA_it}. One can verify that $x^k$ and $\bar{x}^k$ are inflection points for $f_1$, the curvature changing between $\mu_1$ and $L_1$. The inflection points $\bar{x}^k$ result from ensuring continuity of both the function and its gradient values and one can verify that $\bar{x}^k \in [x^k, x^{k+1}]$ for all $k=0,\dots,N-1$.
\end{proof}
The worst-case example from \cref{prop:wc_example_p1} builds on the proof of regime $p_1$ (see \cref{sec:proof_one_step_decrease_dca}). Firstly, condition $G=G^+=L_2 \Delta x$ implies $g_1^k - g_2^k = L_2 (x^k - x^{k+1})$ at each iteration $k=0,\dots,N-1$. Additionally, the worst-case after $N$ iterations implies all gradient norms equal with value of $|U|$ from \cref{prop:wc_example_p1}.

\begin{figure}
    \centering
    \includegraphics[width=0.49\linewidth]{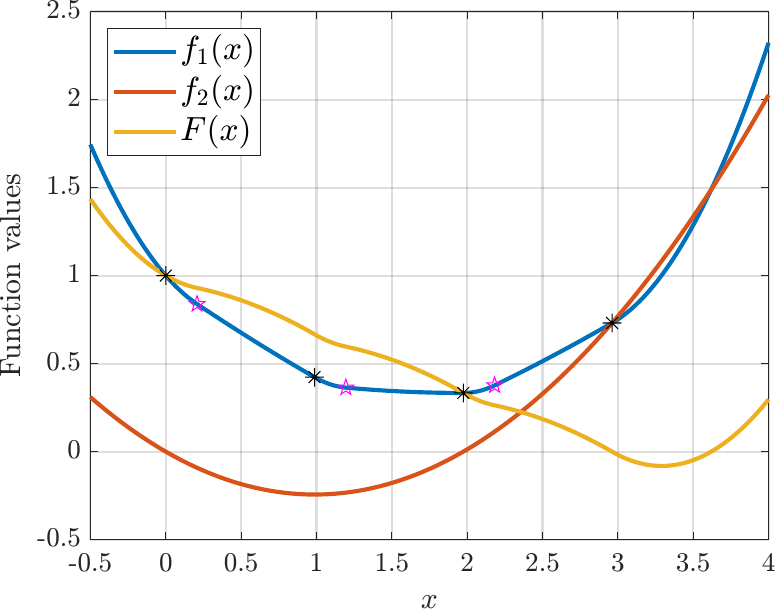}
    \hfill
    \includegraphics[width=0.49\linewidth]{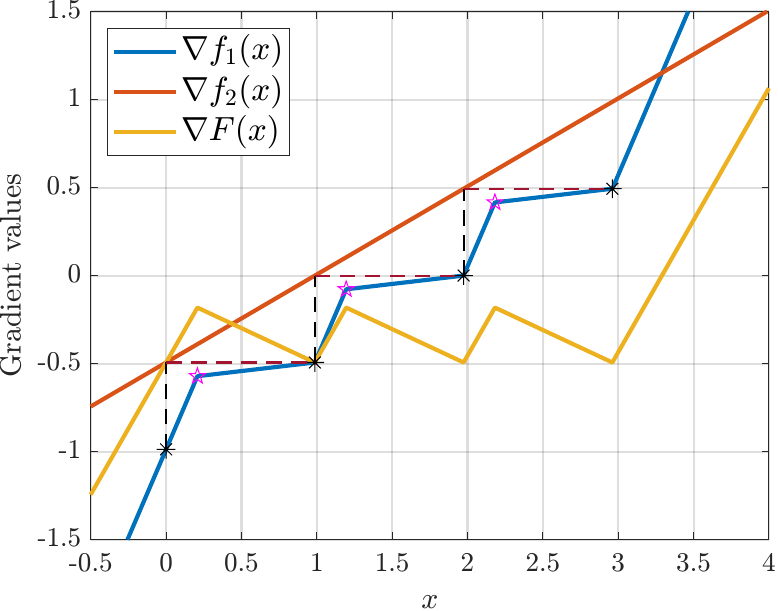}
    \caption{
    Worst-case example for regime $p_1$ after $N=3$ iterations, showing the function values (\textit{left}) and the gradient values (\textit{right}). Setup: $L_1=2$, $\mu_1=0.1$, $L_2=0.5$, $\mu_2=-0.01$. The initial condition is $\Delta = f_1(x^0) - f_2(x^0) = 1$. The iterations begin at $x^0=0$ where $f_1(x^0)=1$ and $f_2(x^0) = 0$. The \textit{black} stars represent the DCA iterations $x^k$, while the \textit{magenta} pentagrams indicate the inflection points of $f_1$, denoted by $\bar{x}^k$, placed between consecutive iterations $x^k$.}
    \label{fig:wc_func_p1}
\end{figure}

\paragraph{DCA iterations.} The plot on the \textit{right} in \cref{fig:wc_func_p1} provides an intuitive illustration of the DCA iteration process. Starting from $x^0=0$, the iteration moves vertically along the y-axis until intersecting the graph of $\nabla f_2$. Next, it moves horizontally along the x-axis until reaching the graph of $\nabla f_1$ at $x^1$. This procedure is then repeated for each subsequent iteration.
 
\subsection{Worst-Case Function Examples for Regime \texorpdfstring{$p_2$}{p2} }
\begin{proposition}[Tightness of $p_2$]\label{prop:wc_example_p2}
Let $f_1^0$, $f_2^0$, $g_2^0$, $x^0$ be some real numbers such that $\Delta \coloneqq f_1^0-f_2^0 > 0$ and let $N\geq 1$ be an integer.
Consider $\mu_1$, $L_1$, $\mu_2$, $L_2$ as belonging to the domain of regime $p_2$ from \cref{tab:DCA_regimes_one_step}, where the conditions are $\mu_1, \mu_2 \geq 0$ and $\max\{\mu_1,\mu_2\} < L_1 < L_2 < \infty$. Let $U{} {}\coloneqq{} {} -\sqrt{\frac{2\Delta}{p_2 N}}$, where the expression of $p_2$ is given in \cref{tab:DCA_regimes_one_step}. Consider the functions $f_1, f_2: \mathbb{R}\rightarrow \mathbb{R}$ defined as:
\begin{align}
    \begin{aligned}
        f_1(x) 
            {} {}\coloneqq{} {}&
        \tfrac{1}{2} L_1 (x-x^0)^2 {}+{}
        g_1^0 (x-x^0) {}+{}
        f_1^0 \\
        f_2(x) {} {}\coloneqq{} {}&
        \left\{
        \def\arraystretch{2}
        \begin{array}{ll} 
        \tfrac{1}{2} \mu_2 (x-x^{0})^2 {}+{} g_2^{0} (x-x^{0}) {}+{}
            f_2^{0} & x \in (-\infty, x^0]; \\
           \tfrac{1}{2} \mu_2 (x-x^k)^2 {}+{}
            g_2^k (x-x^k) {}+{}
            f_2^k  & x \in [x^k, \bar{x}^k];  \\
            \tfrac{1}{2} L_2 (x-x^{k+1})^2 {}+{}
            g_2^{k+1} (x-x^{k+1}) {}+{}
            f_2^{k+1}  & x \in [\bar{x}^k, x^{k+1}]; \\
            \tfrac{1}{2} \mu_2 (x-x^{N})^2 {}+{} g_2^{N} (x-x^{N}) {}+{}
            f_2^{N}    & x \in [x^N, \infty),
        \end{array}            
        \right.
    \end{aligned}
\end{align}%
where for all $k=0,1,\dots,N-1$ it holds $\nabla f_1(x^0)=g_1^0$, $f_2(x^k) = f_2^k$, $\nabla f_2(x^k) = g_2^k$ and:%
\begin{align*}
    \begin{aligned}
        x^{k} {}={}& x^0 - k \frac{U}{L_1}, {}\quad{} \forall k = 0,\dots,N; \\
        g_2^{k} {}={}& g_1^0 - (k+1) U, {}\quad{} \forall k = 0,\dots,N; \\
        f_2^{k} {}={}& f_1(x^{k}) - \frac{N-k}{N} \Delta, {}\quad{} \forall k = 0,\dots,N; \\
        \bar{x}^{k} {}={}& x^k - \frac{L_2 - L_1}{L_2 - \mu_2} \frac{U}{L_1}, {}\quad{} \forall k = 0,\dots,N-1.
    \end{aligned}
\end{align*}%
Then by performing $N$ iterations of \eqref{eq:DCA_it} on the function $F(x)=f_1(x)-f_2(x)$, starting from $x^0$, the following result holds:
\begin{align}
    \tfrac{1}{2} \min_{0\leq k \leq N} \{\|\nabla f_1(x^k) - \nabla f_2(x^k)\|^2\} {}={}
    \frac{\Delta}{p_2(L_1, L_2, \mu_1, \mu_2) N}.
\end{align}
\end{proposition}%
The proof of \cref{prop:wc_example_p2} is similar to the one of \cref{prop:wc_example_p1}.

Regime $p_2$ can be understood as a transformation of regime $p_1$. Specifically, regime $p_2$ appears when the \eqref{eq:DCA_it} iterations are applied in reverse on the function  $-F(x) = f_2(x) - f_1(x)$, with the roles of $f_1$ and $f_2$ interchanged. Thus, the dynamic of regime $p_1$ is mirrored in regime $p_2$.%
\begin{figure}
    \centering
    \includegraphics[width=0.46\linewidth]{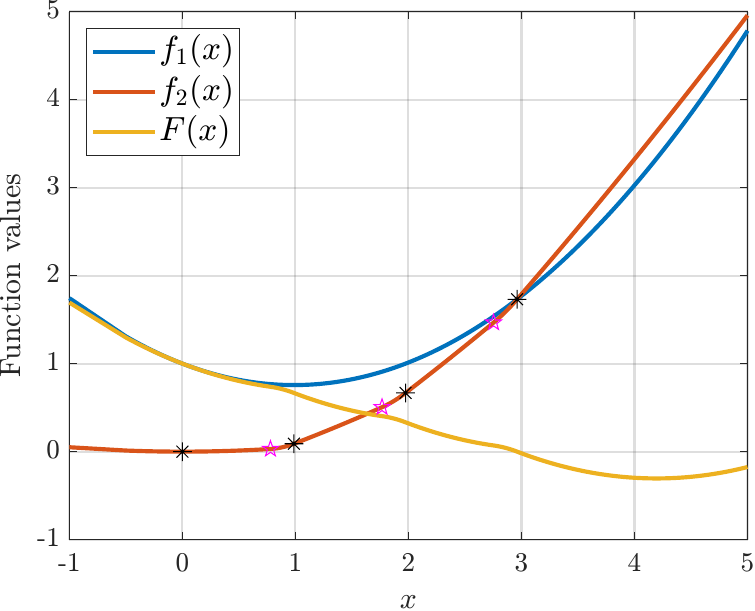}
    \hfill
    \includegraphics[width=0.46\linewidth]{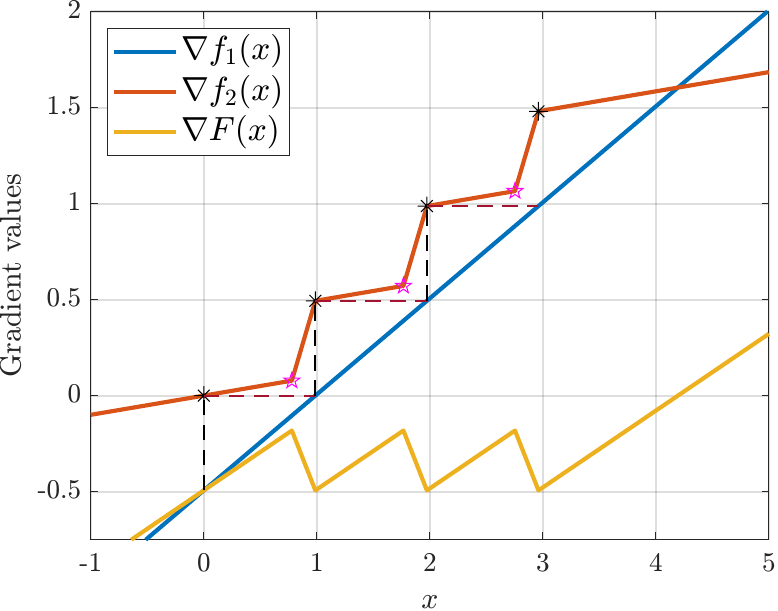}
    \caption{
    Worst-case example for regime $p_2$ after $N=3$ iterations, showing the function values (\textit{left}) and the gradient values (\textit{right}). Setup: $L_1=1.5$, $\mu_1=0.25$, $L_2=2$, $\mu_2=1$. The initial condition is $\Delta = f_1(x^0) - f_2(x^0) = 1$. The iterations begin at $x^0=0$ where $f_1(x^0)=1$ and $f_2(x^0) = 0$. The \textit{black} stars represent the DCA iterations $x^k$, while the \textit{magenta} pentagrams indicate the inflection points of $f_2$, denoted by $\bar{x}^k$, placed between consecutive iterations $x^k$. The dashed lines show the DCA iterations $x^{k+1} = \nabla f_1^*(\nabla f_2(x^k))$.}
    \label{fig:wc_func_p2}
\end{figure}



\clearpage
\section{\MakeUppercase{Decision tree for identifying the active regime}}
\begin{figure}[!ht]
    \centering
\resizebox{\linewidth}{!}{%
\centering
\begin{tikzpicture}[auto]
        \node [box]                                    (p123578)      {$p_1$, $p_2$, $p_3$, $p_4$, $p_5$, $p_6$};
        \node [box, below=4.5cm of p123578, xshift=-3cm]    (p5)    {$p_4$};
        \node [box, right=1cm of p5]                (p7)    {$p_5$};
        \node [box, right=2cm of p7, xshift=-1.1cm]  (p8)    {$p_6$};
        \node [box, right=2cm of p8, xshift=1cm]  (p2)    {$p_2$};
        \node [box, right=2.5cm of p2, xshift=.2cm]   (p3w)    {$p_3$};
        \node [box, right=3cm of p3w, xshift=.5cm]   (p1)    {$p_1$};

        \node [box, below=.5cm of p123578, xshift=4cm]     (p12378)    {$p_1$, $p_2$, $p_3$, $p_5$, $p_6$};
        
        \node [box, below=.5cm of p12378, xshift=4cm]   (p1238)    {$p_1$, $p_2$, $p_3$, $p_6$};
        \node [box, below=.5cm of p1238, xshift=-3cm]  (p238)    {$p_2$, $p_3$, $p_6$};
        \node [box, below=.5cm of p1238, xshift=3cm]   (p13)    {$p_1$, $p_3$};

        \node [box, below=.5cm of p238, xshift=-2.1cm]  (p28)    {$p_2$, $p_6$};
        \draw [->] (p123578) -| (p5) node [midway, above] (TextNode) {$B\leq 0$};
        \draw [->] (p123578) -| (p12378) node [midway, above] (TextNode) {$B > 0$};
                
        \draw [->] (p12378) -| (p7) node [midway, above] (TextNode) {$L_2 \leq \mu_1$ ($F$ s.c.)};
        \draw [->] (p12378) -| (p1238) node [midway, above] (TextNode) {\quad $L_2 > \mu_1$ ($F$ n.c.)};

        \draw [->] (p1238) -| (p238) node [midway, above] (TextNode) {$L_2 > L_1$};
        \draw [->] (p1238) -| (p13) node [midway, above] (TextNode) {$L_1 \geq L_2$};

        \draw [->] (p238) -| (p28) node [midway, above] (TextNode) {$\mu_2 \geq 0$ ($f_2$ s.c.)};
        \draw [->] (p238) -| (p3w) node [pos=0.22, above] (TextNode) {$\mu_2 < 0$ ($f_2$ n.c.)};

        \draw [->] (p28) -| (p8) node [midway, above] (TextNode) {$L_1 \leq \mu_2$ ($F$ ccv.)};
        \draw [->] (p28) -| (p2) node [midway, above] (TextNode) {$L_1 > \mu_2$ ($F$ n.ccv.)};

        \draw [->] (p13) -| (p3w) node [pos=0.2, above] (TextNode) {$E>0$};
        \draw [->] (p13) -| (p1) node [midway, above] (TextNode) {$E \leq 0$};
        
    \end{tikzpicture}
    }
    \caption{Decision tree on selecting regimes after one iteration defined in \cref{thm:one_step_decrease_dca}. Abbreviations -- s.c.: strongly convex; n.c.: nonconvex; ccv.: concave; n.ccv.: nonconcave. Recall that $F=f_1 - f_2$, where we assume $f_1 \in \mathcal{F}_{\mu_1,L_1}$ and $f_2 \in \mathcal{F}_{\mu_2,L_2}$, with $\mu_1 + \mu_2 > 0$ or $\mu_1=\mu_2=0$, and use the notation: $B  {}\coloneqq{}  \mu_1^{-1} + \mu_2^{-1} + L_2^{-1}$ and $E  {}\coloneqq{}  \frac{L_2+\mu_2}{L_1 L_2} \frac{L_2 - L_1}{-\mu_2} + \mu_1^{-1} - L_1^{-1}$.}
    \label{fig:decision_tree_distinguish_6_regimes}
\end{figure}
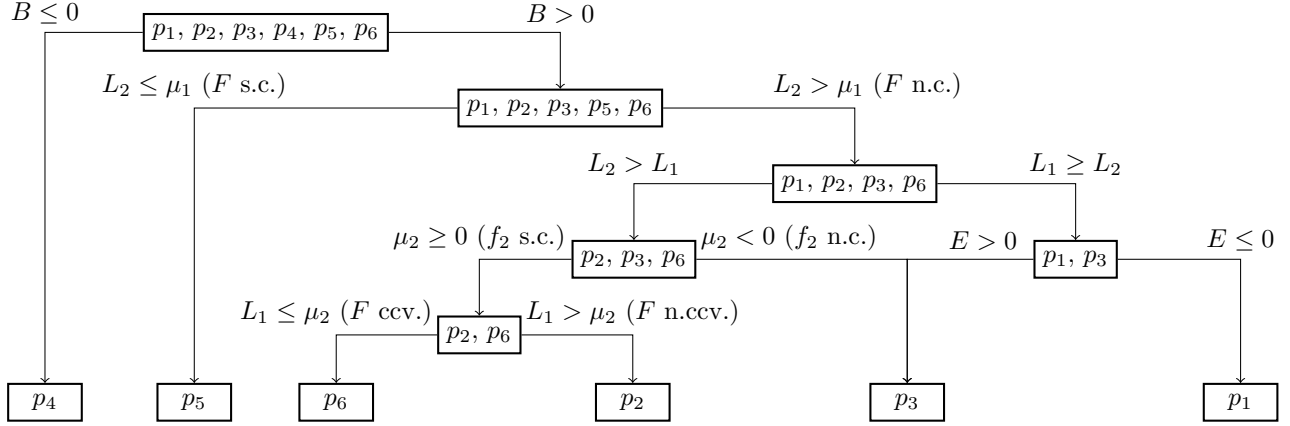

\cref{fig:decision_tree_distinguish_6_regimes} shows a decision tree outlining how to identify the corresponding regime from \cref{thm:one_step_decrease_dca} based on the given curvature parameters.

\section{\MakeUppercase{Sketch of regimes in the nonsmooth cases} (see \texorpdfstring{\cref{corollary:dca_rates_nonsmooth_N_steps}}{Corollary 2})}
\label{app:sec:nonsmooth_plots}

\begin{figure}[!ht]
    \centering
    \includegraphics[width=.4\textwidth]{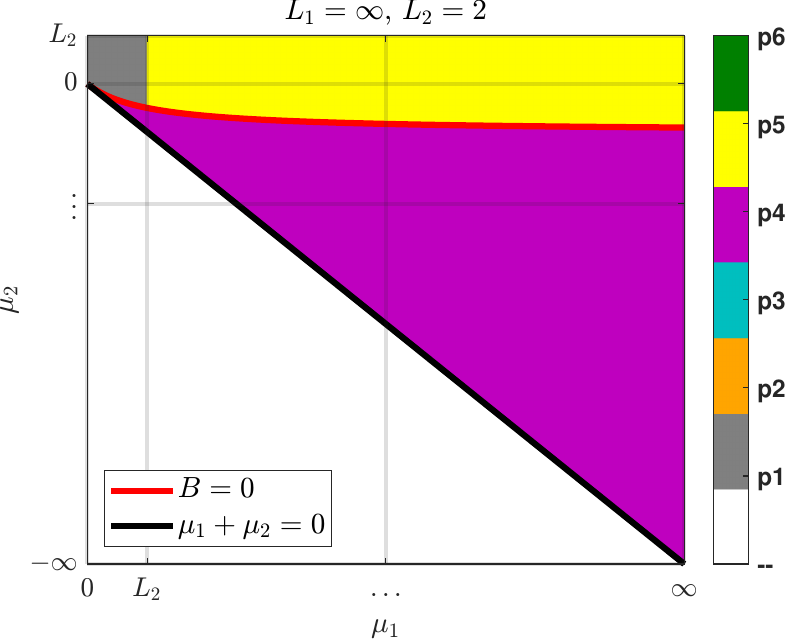}
        {}{}
    \includegraphics[width=.4\textwidth]{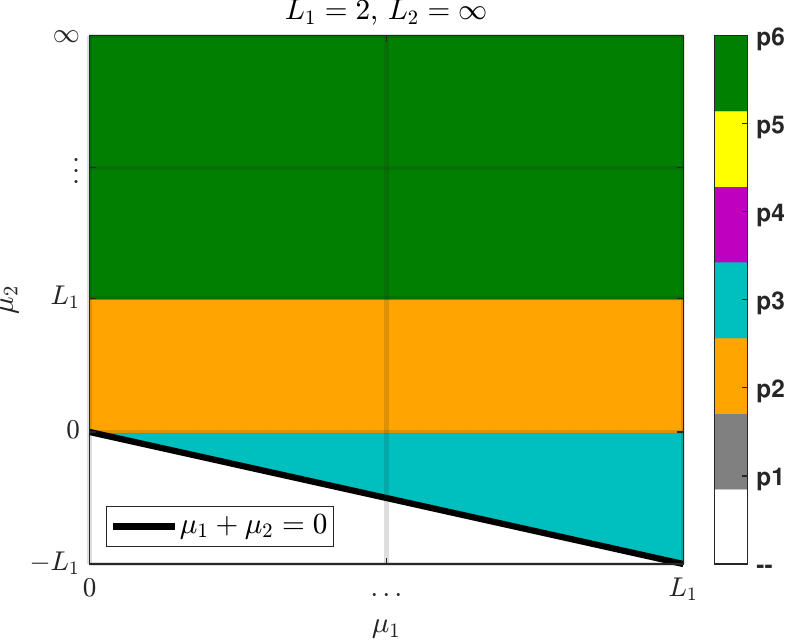}
    \caption{
    All regimes after one \eqref{eq:DCA_it} iteration when either $f_1$ or $f_2$ is nonsmooth (see \cref{tab:DCA_regimes_nonsmooth} for precise expressions). The denominator values $p_i$ are identical for regimes $p_1$ and $p_5$, as well as for regimes $p_2$ and $p_6$. Furthermore, when $L_2 = \infty$ (indicating $F$ (strongly) concave), the threshold $B=\mu_1^{-1} + \mu_2^{-1} + L_2^{-1}=0$ and the limit line $\mu_1 + \mu_2 = 0$ coincide.}
    \label{fig:all_together_one_nonsmooth}
\end{figure}%

\clearpage
\section{\MakeUppercase{On the best curvature splitting problem} (see \texorpdfstring{\cref{subsec:shifting_curvature}}{Section 4})}
In \cref{subsec:shifting_curvature} we examined how to achieve the optimal splitting for a given objective $F=f_1 - f_2$ by subtracting the same curvature term $\lambda \frac{\|\cdot\|^2}{2}$ in both functions. We concluded that if $\mu_1 \leq \mu_2$, the best splitting is achieved by shifting the lower curvature of $f_1$ to $0$, i.e., make it convex. Conversely, when $\mu_1 > \mu_2$, the best splitting is achieved for some $f_2$ weakly convex, with lower curvature $\mu_2 - \lambda < 0$.

\begin{figure}[!ht]
    \centering
    \includegraphics[width=\textwidth]{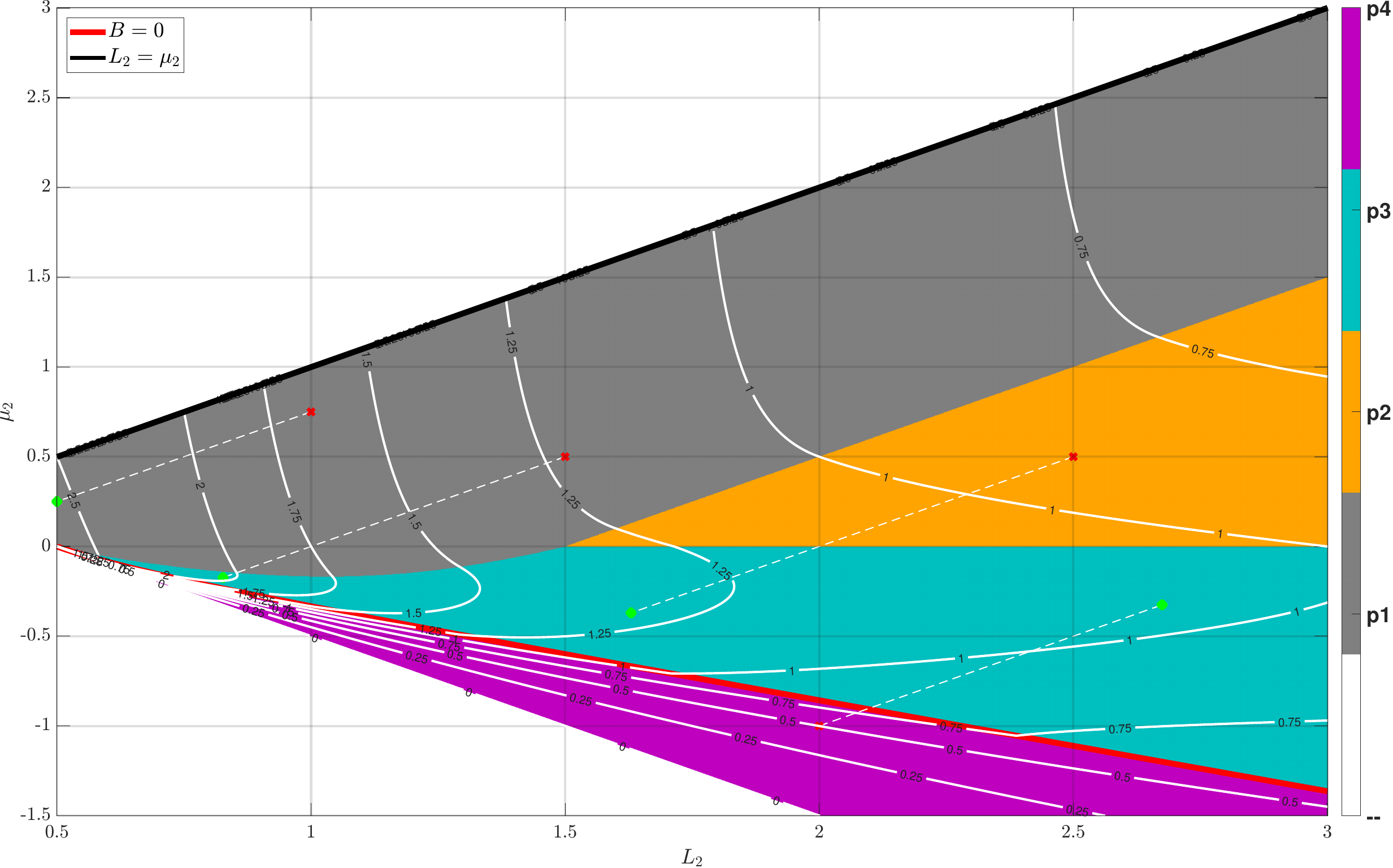}
    \caption{All regimes for a fixed objective function $F$ with $\mu_F = -0.5$ and $L_F=1.5$, along with several mappings of the optimal splittings, shown as transitions from  \textcolor{red}{\textit{red}} to \textcolor{ForestGreen}{\textit{green}} dots along dashed lines with a slope of $1$.}
    \label{fig:case_muF=-0.5_LF=1.5_best_splittings}
\end{figure}

Within this section, we provide additional numerical experiments. In \cref{fig:case_muF=-0.5_LF=1.5_best_splittings} we provide an example with fixed curvature bounds of a nonconvex-nonconcave objective function $F$, namely $\mu_F = -0.5$ and $L_F=1.5$, and illustrate all possible regimes after one iteration. Note that $\mu_F = \mu_1 - L_2$ and $L_F = L_1 - \mu_2$; further on, we examine the regimes based on the ranges of $L_2$ and $\mu_2$. The condition $\mu_1 \geq 0$ implies that $L_2 = \mu_1 - \mu_F \geq -\mu_F$, while the condition $\mu_1 + \mu_2 > 0$ that $\mu_2 > -\mu_1 = L_2 + \mu_F$. The contour lines represent the values of the denominators $p_i$, with $i=1,2,3,4$. The \textcolor{red}{\textit{red}} points mark the initial curvature values of $L_2$ and $\mu_2$, whereas the \textcolor{ForestGreen}{\textit{green}} dots indicate the points with the largest possible $p_i$ obtained through the optimal choice of $\lambda$. Since these shifts are linear in $\lambda$, the dashed lines connecting the dots have a slope of one. 

For example, in the case where $L_2 = 1$ and $\mu_2 = 0.75$, we have $\mu_1 < \mu_2$ and the best splitting is obtained by shifting to the lowest possible value of $L_2$, corresponding to $\mu_1 = 0$. In all other examples, $\mu_1 > \mu_2$ and the optimal splittings are found within regime $p_3$, where $\mu_2 < 0$.

\clearpage
\section{\MakeUppercase{Numerical experiments} (extends \texorpdfstring{\cref{sec:Numerical_experiments}}{Appendix Numerical Experiments})} \label{app:numerical_experiments}

Recall the for any $x \in \mathbb{R}^n$ we define $f_1(x)  {}={}  \kappa \|x\|_1 + \eta \frac{\|x\|^2}{2} + \delta_{\bar{B}(0,1)}(x)$ and $f_2(x)  {}\coloneqq{}  \frac{1}{2}x^T \Sigma x$, 
where $\kappa > 0$, $\eta \geq 0$, $\delta_{\bar{B}(0,1)}$ is the indicator function over the closed Euclidean unit ball, $\Sigma = A^T A$ is the sample covariance matrix, where $A \in \mathbb{R}^{20 n \times n}$, with $n=200$ and $10\%$ sparsity. Note that $\mu_1 = \eta$, $L_1=\infty$, $\mu_2 = \min\{\Lambda(\Sigma)\}$ and $L_2 = \max \{\Lambda(\Sigma)\}$, where $\Lambda(\Sigma)$ denotes the eigenvalues of $\Sigma$. Recall $\tilde{f}_{i}^{\lambda} {}={} f_{i} - \lambda \frac{\|\cdot\|^2}{2}$, with $i=\{1,2\}$, where $\lambda \in \left(-\infty, \eta \right]$ is the curvature shift, being the curvature adjusted functions.


\subsection{Closed-Form Expression of the \texorpdfstring{$\partial \tilde{f}_1^{{\lambda}^{*}}(y)$}{Fenchel conjugate}}
Firstly, for each $y \in \mathbb{R}^n$ we compute $\partial {f_1}^{*}(y) {}={} \argmin_{x \in \mathbb{R}^n} \left\{ \langle x, y \rangle {}-{}  f_1(x) \right\}$. Let $f(x){}\coloneqq{} \kappa \|x\|_1 + \delta_{\bar{B}(0,1)}(x) $, such that $f_1 = f + \eta \frac{\|\cdot\|^2}{2}$ and $\tilde{f}_{1}^{\lambda} {}={} f + (\eta-\lambda) \frac{\|\cdot\|^2}{2}$. From Themelis et al. \citeyearpar[Section V]{Themelis_2020_DC_envelope} we have:
\begin{align*}
    \prox_{\gamma f}(y) {}={} \frac{\sgn(y) \odot \left[|y| - \kappa \gamma \mathbf{1} \right]_{+} }{ \max \left\{1 \,,\, \| \left[|y| - \kappa \gamma \mathbf{1} \right]_{+} \|_2 \right\} },
\end{align*}
where $\gamma>0$, $\left[\cdot\right]_{+}$ denotes the positive part, $\odot$ is the elementwise multiplication, and $\mathbf{1}$ is the vector of all ones. 

\textbf{Case $\eta-\lambda > 0$.} $\tilde{f}_{1}^{\lambda}$ is strongly convex and the subdifferential of the convex conjugate is a singleton. Using the properties of the proximal operator, we get:%
\begin{align*}
    \nabla \tilde{f}_1^{\lambda^*}(y) {}={} \prox_{(\eta-\lambda)^{-1} f}\left((\eta-\lambda)^{-1} y\right).
\end{align*}%
\textbf{Case $\eta - \lambda = 0$.} We get $\tilde{f}_{1}^{\lambda} = f$, which is convex. We obtain
\begin{align*}
    \partial {f}^{*}(y) {}={}
        \left\{
    \def\arraystretch{1.5}
        \begin{array}{ll}
            \mathbf{0} , & \text{ if } y {}={} \mathbf{0} \\
            \mathbf{0} \cup \left\{ x \in \bar{B}(0,1):\, \sgn(x_i) = \sgn(y_i), \, \forall i = 1,\dots,n \right\}, &  \text{ if } |y_i| = \kappa, \, \forall i = 1,\dots,n; \\
            \frac{\sgn(y) \odot [ |y| - \kappa \mathbf{1} ]_{+} }
    { \| \left[ |y| - \kappa \mathbf{1} \right]_{+}  \|_2  }, & \text{ otherwise },
        \end{array}
    \right.
\end{align*}
where $\mathbf{0}$ is the vector of all zeros. Putting together, for all $\lambda \leq \eta$ it holds:%
\begin{align*}
    \partial \tilde{f}_1^{\lambda^*}(y) 
        {}={} 
    \left\{
    \def\arraystretch{1.5}
    \begin{array}{lll}
    \mathbf{0} , & \text{ if } \lambda \leq \eta, & y {}={} \mathbf{0}; \\
    \mathbf{0} \cup \left\{ x \in \bar{B}(0,1): \sgn(x_i) = \sgn(y_i), \, \forall i = 1,\dots,n \right\}, &  \text{ if } \lambda = \eta, & |y_i| = \kappa, \, \forall i = 1,\dots,n; \\
        \frac{\sgn(y) \odot [ |y| - \kappa \mathbf{1} ]_{+} }
    {\max \{ \eta - \lambda , \| [ |y| - \kappa \mathbf{1} ]_{+} \|_2 \} }, &  \text{ otherwise.}
    \end{array}
    \right.
\end{align*}
The subdifferential is a singleton unless $\lambda = \eta$ and $y$ is on the boundary of the $l_{\infty}$-ball of radius $\kappa$. In this case, in the numerical experiments we select the subgradient $\mathbf{0} \in \partial \tilde{f}_1^{\lambda^*}(y)$.

\subsection{Extended Experiment Results}
Let $\bar{M} \leq M = 1000$ be the number of runs converging to the same solution (non-trivial and with the desired sparsity) and denote by $N_{\varepsilon}$ the average number of iterations to reach a certain accuracy $\varepsilon = \{10^{-1}, 10^{-2}, \dots, 10^{-12}\}$, defined as
    $
        N_{\varepsilon}  {}\coloneqq{}  \frac{1}{\bar{M}} \sum_{d=1}^{\bar{M}} \argmin_{0 \leq k \leq N} \{  \|\tilde{g}_1^k - \nabla \tilde{f}_2^{\lambda}(x^k)\|^2 \leq \varepsilon \}.
    $

\textbf{Case 1}: $\eta > \mu_2$. Parameters: $\eta =\mu_1 = 0.5$, $\kappa = 0.02$ (400 runs kept out of 1000). Here, $\lambda^* = 0.4413$ and the maximum \textit{theoretical} shift is $\lambda_{\text{max}}=\frac{\mu_1+\mu_2}{2} = 0.4441$. With $\lambda > \mu_2 = 0.3882$, $\tilde{f}^{\lambda}_2$ becomes weakly convex. The results are reported in \cref{tab:SPCA_experiment_eta_gg_mu2}. Using $\lambda^*$ achieves at least a twofold acceleration compared to $\lambda=0$. Higher $\lambda$ values, as $\lambda_{\text{max}}$, making $\tilde{f}^{\lambda}_2$ even more nonconvex, further improve the algorithm. Conversely, adding curvature to both functions ($\lambda<0$) slows the convergence.%

\begin{table}[!ht]
\centering
\caption{\textbf{(Case 1)} Average number of iterations $N_{\varepsilon}$ to reach $\varepsilon$ accuracy of the minimum subgradient norm in the setting $\eta=\mu_1=0.5$, $\kappa=0.02$, corresponding to $\mu_1 > \mu_2$. Lowest values are obtained for $\textcolor{blue}{\lambda = \lambda_{\text{max}}=0.4441}$, whereas $\textcolor{red}{\lambda^* = 0.4413}$ performs the second best.}
\label{tab:SPCA_experiment_eta_gg_mu2}
\resizebox{\linewidth}{!}{%
    \begin{tabular}{@{}r|cccccccccccc}
 \backslashbox[0pt][l]{ $\lambda$}{$\varepsilon$}
 & $10^{-1}$  & $10^{-2}$ &  $10^{-3}$ & $10^{-4}$ &  $10^{-5}$ &  $10^{-6}$ &  $10^{-7}$ &  $10^{-8}$ &  $10^{-9}$ &  $10^{-10}$ &  $10^{-11}$ & $10^{-12}$ \\ \hline
 0       & 1.00 & 5.74 & 17.73 & 46.70 & 96.70  & 150.17 & 242.57 & 336.27 & 439.46 & 544.63 & 651.33  & 758.74  \\
\textcolor{red}{$\lambda^* {}={}$ 0.4413}  & \textcolor{red}{1.00} & \textcolor{red}{3.65} & \textcolor{red}{8.44}  & \textcolor{red}{21.49} & \textcolor{red}{44.76}  & \textcolor{red}{70.75}  & \textcolor{red}{112.57} & \textcolor{red}{154.93} & \textcolor{red}{200.98} & \textcolor{red}{247.93} & \textcolor{red}{295.43}  & \textcolor{red}{343.06}  \\
$-0.5\lambda^* {}={}$ -0.2207 & 1.00 & 6.83 & 22.34 & 59.32 & 122.73 & 189.31 & 307.17 & 426.39 & 558.51 & 693.40 & 830.17  & 967.74  \\
$-\lambda^*{}={}$ -0.4413 & 1.00 & 8.05 & 27.00 & 72.30 & 147.77 & 229.13 & 375.34 & 521.24 & 681.09 & 844.69 & 1010.64 & 1177.50 \\
$0.5\lambda^*{}={}$ 0.2207  & 1.00 & 4.54 & 13.06 & 33.91 & 70.19  & 108.43 & 176.35 & 244.05 & 318.63 & 394.92 & 472.42  & 550.15  \\
\textcolor{blue}{$\lambda_{\text{max}} {}={}$ 0.4441}  & \textcolor{blue}{1.00} & \textcolor{blue}{3.63} & \textcolor{blue}{8.35}  & \textcolor{blue}{21.23} & \textcolor{blue}{44.01}  & \textcolor{blue}{70.38}  & \textcolor{blue}{111.91} & \textcolor{blue}{153.78} & \textcolor{blue}{199.95} & \textcolor{blue}{246.61} & \textcolor{blue}{293.76}  & \textcolor{blue}{340.99} \\  \hline
    \end{tabular}}
\end{table}

\textbf{Case 2}: $\eta < \mu_2$. Parameters: $\eta=0.2$, $\kappa=0.02$ (703 runs kept out of 1000). We get $\lambda^* = \mu_1 = 0.2$. Using $\lambda^*$ to make $\tilde{f}^{\lambda}_1$ convex improves the convergence speed by approximately $20\%$ compared to the initial splitting. Generally, decreasing the curvature of $f_2$ improves the convergence.

\begin{table}[!ht]
\centering
\caption{\textbf{(Case 2)} Average number of iterations $N_{\varepsilon}$ to reach $\varepsilon$ accuracy of the minimum subgradient norm in the setting $\eta=0.2$, $\kappa=0.02$, corresponding to $\eta < \mu_2$. Lowest values are obtained for \textcolor{red}{$\lambda^*= \mu_1 = 0.2$}.\\}
\label{tab:SPCA_experiment_lambda_opt_eq_mu1}
\resizebox{\linewidth}{!}{%
    \begin{tabular}{@{}r|cccccccccccc}
 \backslashbox[0pt][l]{ $\lambda$}{$\varepsilon$}
 & $10^{-1}$  & $10^{-2}$ &  $10^{-3}$ & $10^{-4}$ &  $10^{-5}$ &  $10^{-6}$ &  $10^{-7}$ &  $10^{-8}$ &  $10^{-9}$ &  $10^{-10}$ &  $10^{-11}$ & $10^{-12}$ \\ \hline
0 & 1.04 & 5.88	& 19.21 & 49.85 & 99.16 & 153.46 & 245.40 & 339.67 & 442.23 & 547.61 & 654.99 & 763.01 \\
\textcolor{red}{$\lambda^*$ = 0.2} & \textcolor{red}{1.04} & \textcolor{red}{4.32} & \textcolor{red}{14.23} & \textcolor{red}{37.10} & \textcolor{red}{75.09} & \textcolor{red}{115.90} & \textcolor{red}{185.48} & \textcolor{red}{256.33} & \textcolor{red}{332.81} & \textcolor{red}{412.39} & \textcolor{red}{492.94} & \textcolor{red}{574.18} \\
$-0.5\lambda^*$ = -0.1 & 1.04 & 6.48 & 21.62 & 55.52 & 110.98 & 172.24 & 276.85 & 382.11 & 498.04 & 616.29 & 737.17 & 859.14 \\ 
$-\lambda^*$ = -0.2 & 1.04 & 7.06 & 23.87 & 61.44 & 123.41 & 192.46 & 307.80 & 426.00 & 554.21 & 686.04 & 820.61 & 956.00 \\
$0.5\lambda^*$ = 0.1 & 1.04 & 5.26 & 16.83 & 43.47 & 87.38 & 135.39 & 215.03 & 298.14 & 387.83 & 479.93 & 574.02 & 668.76  \\ \hline
    \end{tabular}}
\end{table}


\begin{figure}[!ht]
    \centering
    \includegraphics[width=.7\linewidth]{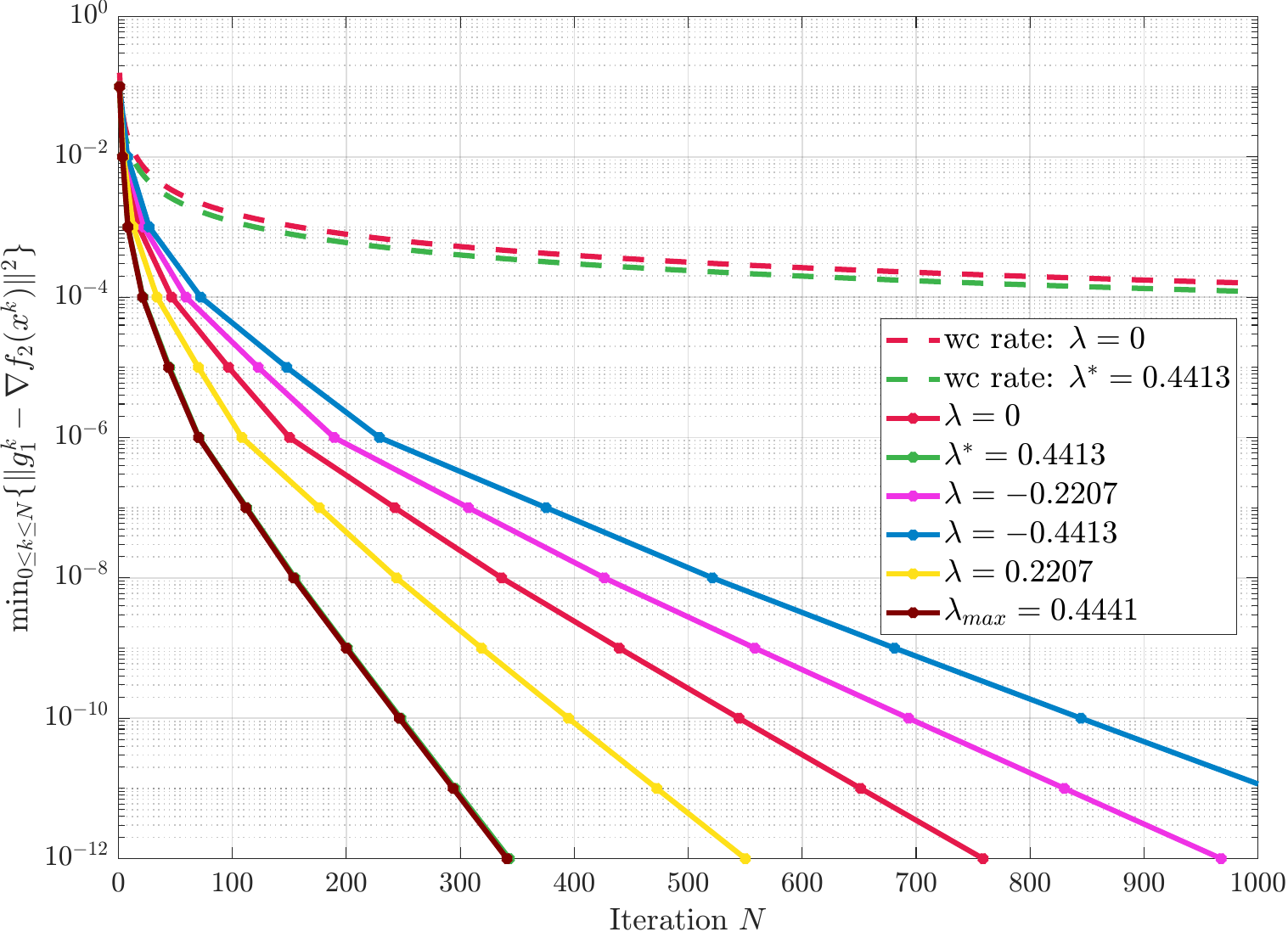}
    \caption{\textbf{Case 1.} Setting: $\eta=0.5 > \mu_2$, $\kappa=0.02$, $\lambda = \{0, \pm \lambda^*, \pm 0.5 \lambda^*, \lambda_{\text{max}}\}$. Lowest values are obtained for $\lambda_{\max}=0.4441$. Average number of iterations $N_{\varepsilon}$ to reach $\varepsilon$ accuracy of the best subgradient norm $\|g_1^k - \nabla f_2(x^k)\|$, where $g_1^k \in \partial f_1(x^k)$. Although the actual rate exceeds theoretical guarantees (\textit{wc rate}), there is correlation between optimizing the worst-case for $\lambda$ and better performance in the experiment.}
    \label{fig:SPCA_experiment_eta_gg_mu2}
\end{figure}

\begin{figure}[!ht]
    \centering
    \includegraphics[width=.7\linewidth]{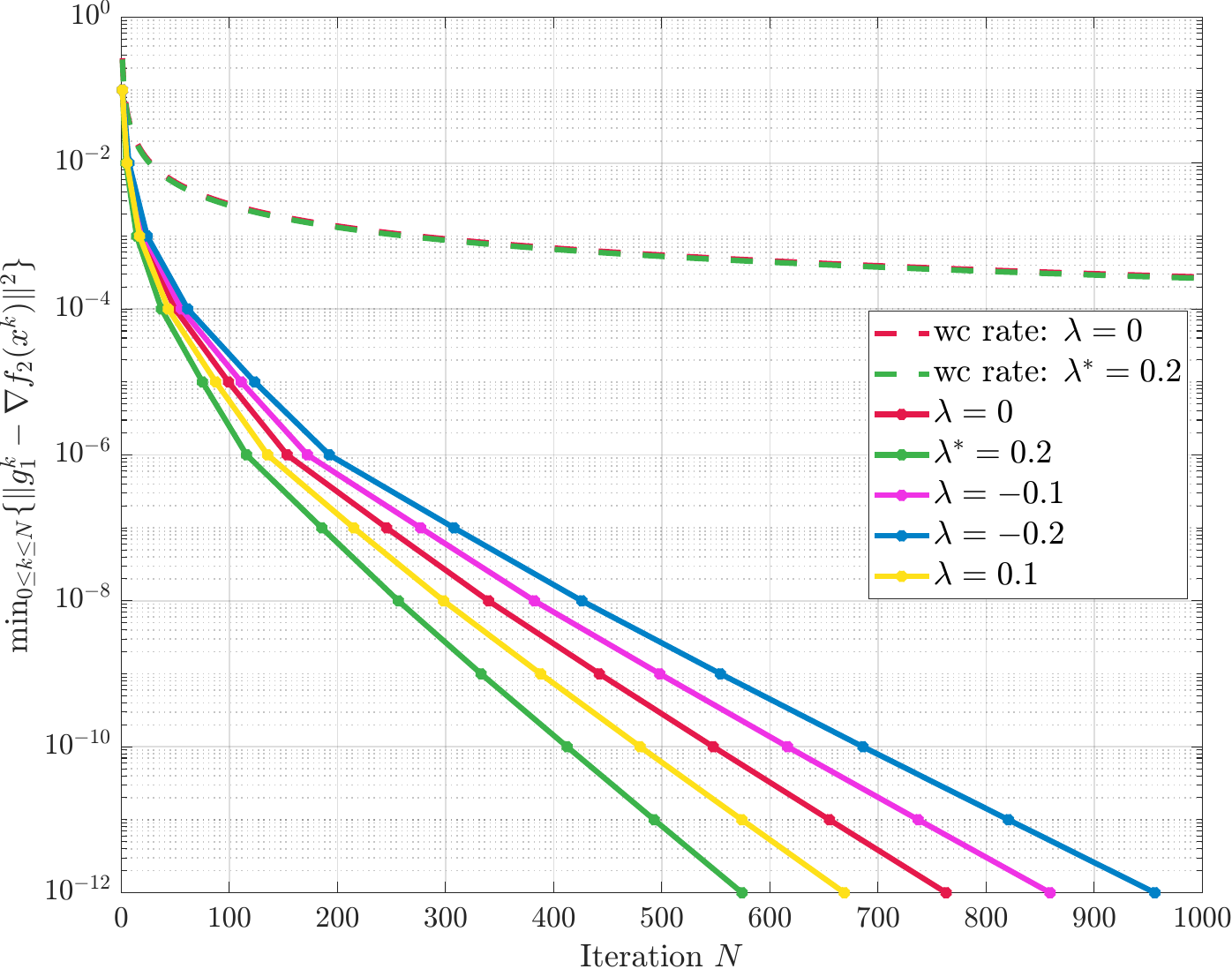}
    \caption{\textbf{Case 2.} Setting: $\eta=0.2 = \mu_1 < \mu_2$, $\kappa=0.02$, $\lambda = \{0, \pm \lambda^*, \pm 0.5 \lambda^*\}$. Average number of iterations $N_{\varepsilon}$ to reach $\varepsilon$ accuracy of the best subgradient norm $\|g_1^k - \nabla f_2(x^k)\|$, where $g_1^k \in \partial f_1(x^k)$.  Lowest values are obtained for $\lambda^* = \lambda_{\text{max}} = \mu_1 = 0.2$. There is correlation between optimizing the worst-case for $\lambda$ and better performance in the experiment.} 
    \label{fig:SPCA_experiment_lambda_opt_eq_mu1}
\end{figure}%
\section{\MakeUppercase{Github repository}}
We provide in the \href{https://github.com/teo2605/DCA_AISTATS25}{GitHub repository} 
\textsc{Matlab} scripts to support the numerical conjectures and to reproduce all the simulations, figures and tables presented in the paper.
\end{document}